\documentclass[]{article}

\usepackage[margin=1in]{geometry}

\usepackage{bbm}
\usepackage{xcolor}
\usepackage{mathrsfs}
\usepackage{tikz}
    \usetikzlibrary{arrows}    
\usepackage{float}
\usepackage{algorithm}
\usepackage[noend]{algpseudocode}
\usepackage{comment}
\usepackage{amsmath}
\usepackage{bm}
\usepackage{amssymb}
\usepackage{amsthm}
\usepackage{color}
\usepackage{graphicx}
\usepackage{multicol,multirow}
\usepackage{subfigure}
\usepackage{todonotes}

\usepackage{natbib}
 \bibpunct[, ]{(}{)}{,}{a}{}{,}%
\usepackage{hyperref}\hypersetup{
    colorlinks,
    linkcolor={red!50!black},
    citecolor={blue!50!black},
    urlcolor={blue!80!black}
}

\usepackage{amsmath,amssymb,enumitem,xspace, calc}

\newcommand{\Z}{\mathbb{Z}}

\newcommand{\R}{\mathbb{R}}

\newcommand{\numberthis}{\addtocounter{equation}{1}\tag{\theequation}}

\DeclareMathOperator{\conv}{conv}

\DeclareMathOperator{\cl}{cl}

\newcommand{\norm}[2]{\left\|#2\right\|_{#1}}

\usepackage{blkarray}

\newcommand{\ignore}[1]{{ }}

\usepackage{tikz}
\usetikzlibrary{positioning}
\usetikzlibrary{calc}
\usetikzlibrary{arrows.meta}
\usetikzlibrary{decorations.markings}
\usetikzlibrary{shapes.misc}
\usetikzlibrary{matrix,shapes,arrows,fit,tikzmark}

\tikzset{   
	every picture/.style={remember picture,baseline},
	every node/.style={anchor=base,align=center,outer sep=0.pt},
	every path/.style={thick},
}

\newcommand{\dd}{\mathcal{D}}
\newcommand{\Nodes}{N}
\newcommand{\A}{A}
\newcommand{\rootnode}{r}
\newcommand{\terminalnode}{t}
\newcommand{\head}[1]{h_{#1}}
\newcommand{\tail}[1]{t_{#1}}

\newcommand{\feasibleSet}{Z}
\newcommand{\val}[1]{\nu_{#1}}
\newcommand{\len}[1]{l_{#1}}
\newcommand{\valvec}{\bm{\val{}}}
\newcommand{\lenvec}{\bm{\len{}}}
\newcommand{\U}[1]{u_{#1}}
\newcommand{\Uvec}{\bm{\U{}}}

\DeclareMathOperator*{\cond}{cond}
\newcommand{\QDiagMax}{\ensuremath{Q_{\max}}}
\newcommand{\VDiagMax}{\ensuremath{V_{\max}}}

\newtheorem{theorem}{Theorem}
\newtheorem{proposition}{Proposition}
\newtheorem{corollary}{Corollary}
\theoremstyle{definition}
\newtheorem{definition}{Definition}
\newtheorem{lemma}{Lemma}
\theoremstyle{remark}
\newtheorem{remark}{Remark}
\newtheorem{example}{Example}

\newcommand*{\defeq}{\stackrel{\text{def}}{=}}

\tikzstyle{zero arc} = [draw,dashed, line width=0.5pt,->]
\tikzstyle{one arc} = [draw,line width=0.5pt,->]
\tikzstyle{zero arc infeasible} = [zero arc,
preaction={
	draw,gray!30!white,-,
	double=gray!30!white,
	double distance=3pt,
}]
\tikzstyle{one arc infeasible} = [one arc,
preaction={
	draw,gray!30!white,-,
	double=gray!30!white,
	double distance=3pt,
}]
\tikzstyle{main node} = [circle,fill=gray!60,font=\scriptsize, inner sep=1pt, text=gray!60]
\tikzstyle{text node} = [font=\scriptsize]
\tikzstyle{arc text} = [font=\scriptsize]
\tikzstyle{optimal arc} = [-,line width=3pt, color=blue!20!white]
\tikzstyle{blocked arc} = [-,line width=3pt, color=red]

\begin{document}

\title{Real-time solution of quadratic optimization problems with banded matrices and indicator variables}

\author{Andr\'es G\'omez\thanks{Industrial \& Systems Engineering, University of Southern California, gomezand@usc.edu}, Shaoning Han\thanks{Industrial \& Systems Engineering, University of Southern California, shaoning@usc.edu} and Leonardo Lozano\thanks{Operations, Business Analytics \& Information Systems, University of Cincinnati, leolozano@uc.edu}}

\date{May 2024}


\maketitle

\begin{abstract}

We consider mixed-integer quadratic optimization problems with banded matrices and indicator variables. These problems arise pervasively in statistical inference problems with time-series data, where the banded matrix captures the temporal relationship of the underlying process. In particular, the problem studied arises in monitoring problems, where the decision-maker wants to detect changes or anomalies. We propose to solve these problems using decision diagrams. In particular we show how to exploit the temporal dependencies to construct diagrams with size polynomial in the number of decision variables. We also describe how to construct the convex hull of the set under study from the decision diagrams, and how to deploy the method online to solve the problems in milliseconds via a shortest path algorithm. 
\end{abstract}

\section{Introduction}

We consider mixed-integer quadratic optimization (MIQO) problems of the form 
\begin{subequations}\label{eq:miqo}
\begin{align}
\min_{\bm{x}\in \R^n, \bm{z}\in \{0,1\}^n, x_0\in \R}\;&\bm{d^\top x}+\bm{c^\top z}+\frac{1}{2}x_0\label{eq:miqo_objective}\\
\text{s.t.}\;&\bm{x^\top Qx} \leq x_0\label{eq:miqo_epigraph}\\
& \bm{x}\circ(\bm{1}-\bm{z})=\bm{0}\label{eq:miqo_indicators}\\
&\bm{z}\in Z,\label{eq:miqo_constraints}
\end{align}
\end{subequations}
where $\bm{x}$ are the continuous decision variables, $\bm{z}$ are the indicator variables, $x_0$ is an auxliary variable representing the epigraph of the nonlinear term, matrix $\bm{Q}\in \R^{n\times n}$ is positive definite, ``$\circ$" denotes the Hadamard (entrywise) product, $\bm{0}$ and $\bm{1}$ are $n$-dimensional vectors of zeros and ones, indicator constraints \eqref{eq:miqo_indicators} impose the logical relationship $x_i\neq 0\implies z_i=1$, and $Z\subseteq \{0,1\}^{n}$ encodes additional constraints on the indicator variables. Even if constraints \eqref{eq:miqo_constraints} are removed, problem \eqref{eq:miqo_objective}-\eqref{eq:miqo_indicators} is challenging to solve due to the presence of the non-convex indicator constraints. In this paper we focus on the case where $\bm{Q}$ is a \emph{banded matrix with bandwidth $k\in \Z_+$}, that is, $Q_{ij}=0$ if $|i-j|>k$. MIQO problems with banded matrices often arise in inference problems with time-series data. For example, slowly varying sparse regression problems \citep{bertsimas2024slowly}, sparse and smooth denoising problems \citep{atamturk2021sparse}, and outlier detection problems in time series \citep{gomez2021outlier}, all fall under the umbrella of quadratic optimization with banded matrices.

Our main motivation to study problem \eqref{eq:miqo} with banded matrices is to tackle inference problems arising from the real-time monitoring of data streams, as described next. In such settings, problem \eqref{eq:miqo} needs to be solved in \emph{real-time}, precluding the use of standard mixed-integer optimization solvers. In fact, problem \eqref{eq:miqo} is rarely solved in such settings, with researchers and practitioners typically resorting to heuristics or simple convex approximations instead.

\subsection{Monitoring problems} 
Monitoring systems are increasingly prevalent in manufacturing systems \citep{yan2017anomaly} and personalized medicine \citep{dunn2018wearables}, due to ease of access to sensors and wearables. Processing this data, inferring the values of the process under study, and detecting anomalous behavior calls for the solution of optimization problems of the form 
\begin{align}\label{eq:inference}
\min_{\bm{x}\in \R^n}\;& \frac{1}{2}\sum_{i=1}^n (y_i-x_i)^2 +\frac{1}{2}\bm{x^\top Rx}+\mu \|\bm{x}\|_0,
\end{align}
where $\{y_i\}_{i=1}^n$ is the time series data generated by the sensor, $\bm{R}\succeq \bm{0}$ is a suitable regularization (banded) matrix capturing the temporal evolution of the process, $\mu>0$ is a regularization parameter to be tuned and $\|\bm{x}\|_0$ is the so-called $\ell_0$-norm, corresponding to the support of vector $\bm{x}$. Problem \eqref{eq:inference} is a special case of problem \eqref{eq:miqo} with $\bm{c}=\bm{1}\mu$, $\bm{d}=-\bm{y}$, $\bm{Q}=\bm{I}+\bm{R}$, $Z=\{0,1\}^n$, a constant term $\frac{1}{2}\|\bm{y}\|_2^2$ added to the objective function, and $\|\bm{x}\|_0$ is represented as the summation of the indicator variables after adding constraint \eqref{eq:miqo_indicators}. Efficiently solving \eqref{eq:inference} offline is already a non-trivial task, since (despite substantial recent improvements) technology for MIQO and mixed-integer nonlinear optimization in general is lagging behind technology for mixed-integer linear optimization (MILO). In monitoring systems, instances of problem \eqref{eq:inference} need to be solved online, each time a new batch of data is generated by the sensor. The resulting optimization problem could include the complete time series $\bm{y}$ generated up to that point, or only the most recent $n$ observations. In either case, with modern sensors generating new information in seconds or even milliseconds, the time frame to solve \eqref{eq:inference} is short indeed.

We point out that problem \eqref{eq:inference} is rarely encountered explicitly in the literature, due in part to the perceived intractability of solving it in real-time. However, variants involving  $\ell_1$-approximations of the $\ell_0$-norm are pervasive in the quality control, statistical and machine learning literatures \citep[e.g., see][]{chen2001atomic,yan2017anomaly,friedrich2017fast,vogelstein2010fast,lin2014alternating,candes2008enhancing,donoho2005stable,tibshirani2005sparsity,rinaldo2009properties,kim2009ell_1,mammen1997locally,ruppert2002selecting,yan2014image,zou2009multivariate,guo2016spline}. Indeed, after replacing the non-convex $\|\bm{x}\|_0$ term with the norm $\|\bm{x}\|_1$, problem \eqref{eq:inference} is convex and can be solved easily. Unfortunately, such approximations may also lead to subpar estimators \citep{bertsimas2016best,mazumder2023subset}.

We now review common regularization terms $\bm{x^\top Rx}$. A common approach is to use multiples of the $k$-th order differences: given $i\in \Z_+$, let $\Delta_i^{(0)}\bm{x}=x_{i}$ be the zeroth-order difference, and let the $k$th-order difference be defined recursively as $\Delta_i^{(k)}\bm{x}=\Delta_{i+1}^{(k-1)}\bm{x}-\Delta_i^{(k-1)}\bm{x}$. Then we can set the regularization term as \begin{equation}\label{eq:kth-Diff}\bm{x^\top Rx}=\lambda\sum_{i=1}^{n-k}\left(\Delta_i^{(k)}\bm{x}\right)^2,\end{equation} where $\lambda\geq 0$ is a smoothness parameter to be tuned. For example, if $k=0$, then $\bm{x^\top R x}=\lambda \|\bm{x}\|_2^2$ is the ridge regularization \citep{hoerl1970ridge} and matrix $\bm{R}$ is diagonal; if $k=1$, then $\bm{x^\top R x}=\lambda \sum_{i=1}^{n-1} (x_{i+1}-x_i)^2$ enforces the smoothness term used by \cite{atamturk2021sparse}, and the matrix $\bm{R}$ is tridiagonal; if $k=2$, then $\bm{x^\top R x}=\lambda \sum_{i=1}^{n-2} (x_{i+2}-2x_{i+1}+x_i)^2$ is the Hodrick-Prescott filter \citep{hodrick1997postwar}, and the matrix has bandwidth 2. More generally, using the $k$th order differences yields a matrix of bandwidth $k$. Another common regularization term commonly used is simply a moving average filter given by 
\begin{equation}\label{eq:movingAverage}\bm{x^\top Rx}=\lambda\sum_{i=1}^{n}\left(x_i-\frac{1}{\min\{k,i-1\}}\sum_{j=1}^{\min\{k,i-1\}} x_{i-j}\right)^2\end{equation}
for some prespecified parameter $k\in \Z_+$, also yielding a matrix of bandwidth $k$. Other regularization terms used include simple modifications of the above, e.g., including coefficients in the moving average filter to give more weight to recent observations.

Approaches relying on $\ell_1$ approximations do not use binary variables and thus cannot impose additional logical constraints \eqref{eq:miqo_constraints}, but such constraints can be invaluable to include a variety of priors in the inference problem.
 The most common prior includes hard constraints on the number of non-zero values, $\sum_{i=1}^n z_i\leq \kappa$ (in which case typically $\bm{c}=\bm{0}$). However, more sophisticated priors such as forcing non-zero values to come in batches of at least $\tau$ consecutive periods can also be enforced. 

Problem \eqref{eq:inference} with zeroth-order differences in \eqref{eq:kth-Diff} is separable and can be solved easily. Moreover, it can be reformulated as a continuous second-order cone program (SOCP) using the perspective reformulation \citep{frangioni2006perspective,gunluk2010perspective}, an approach that is now standard in mixed-integer optimization \citep{hazimeh2022sparse,xie2020scalable} -- in fact, commercial off-the-shelf solvers now perform these reformulations automatically to improve the convex relaxations of MIQO problems.
 Problem \eqref{eq:inference} with the first order differences in \eqref{eq:kth-Diff}, and thus tridiagonal matrices, has also received detailed attention in the literature.  \citet{atamturk2021sparse} propose strong formulations and decomposition algorithms for the problem, exploiting the substructures involving two variables only. \citet{liu2023graph} show that, without additional constraints \eqref{eq:miqo_constraints}, problem \eqref{eq:miqo} with tridiagonal matrices is polynomial-time solvable via a dynamic programming algorithm. In related works, \citet{jewell2018exact} show that a dynamic programming algorithm solves a similar problem to \eqref{eq:inference} with constraints and nonlinear terms coupling consecutive variables only. Alternatively, problem \eqref{eq:inference} with the first function class can be interpreted as an mixed-integer quadratic optimization problem with a Stieltjes matrix, which is polynomial time solvable without additional constraints \citep{han2022polynomial}. Finally, since the sparsity pattern of the matrix induces a tree, the problem is also polynomial time solvable if there are no additional constraints \citep{bhathena2024parametric}, or with a cardinality constraint but provided that all coefficients $c_i$ are identical \citep{das2008algorithms}.

 All the aforementioned properties on the matrix (tridiagonal, tree, Stieltjes) disappear for bandwidths greater than one. While strong formulations have been proposed for problems with bandwidth $k\geq 2$ that improve upon perspective-based relaxations when used with off-the-shelf solvers \citep{han2021compact}, the formulations are neither ideal nor do they lead to real-time solution approaches for problems with $n$ in the hundreds.  
 
 \subsection{Overview of the proposed method} 
 Our proposed approach to tackle problem \eqref{eq:inference} is based on two fundamental ideas. The first idea relates to the convexification of the problem. Letting $\bm{Q}=\bm{I}+\bm{R}$ and defining $$X_Q\defeq\left\{\bm{x}\in \R,\bm{z}\in \{0,1\}^n,x_0\in \R: \bm{x^\top Qx}\leq x_0,\; \bm{x}\circ (\bm{1}-\bm{z})=\bm{0}\right\},$$
 it follows from standard arguments that \eqref{eq:inference} is equivalent to the convex optimization problem
 \begin{align}\label{eq:inferenceConvex}
\frac{1}{2}\|\bm{y}\|_2^2+\min_{\bm{x},\bm{z},x_0}\;&-\bm{y^\top x}+\mu\left(\bm{1^\top z}\right)+\frac{1}{2}x_0\;
\text{     s.t.  }(\bm{x},\bm{z},t)\in \text{conv} (X_Q),
 \end{align}
 where ``$\conv$" denotes the convex hull of a set. 
 In itself, the convex representation \eqref{eq:inferenceConvex} is not helpful, as it requires the computation of the convex hull of a non-polyhedral set, a task that is as difficult (if not more) than solving \eqref{eq:inference} in the first place. However, the critical observation is that \emph{set $X_Q$ does not depend on the data $\bm{y}$}. Therefore, in principle, the convexification can be done offline, and then the real-time solution of \eqref{eq:inference} reduces to a convex optimization problem, which can be performed faster. In this paper, we provide an algorithm to compute an SOCP representation of $\text{cl conv} (X_Q)$ in an extended formulation when $\bm{Q}$ is a banded matrix. Interestingly, the runtime of this algorithm is typically faster, in some cases by orders-of-magnitude, than the solution of a single problem \eqref{eq:inference} via off-the-shelf solvers. Equipped with this representation, the same off-the-shelf solver requires only a few seconds to solve \eqref{eq:inferenceConvex}, while requiring several minutes to solve \eqref{eq:inference} from scratch.

\begin{figure}[!h]
    \begin{center}
    \subfigure[$k=2$]{
    \includegraphics[width=0.45\textwidth,trim={11cm 6cm 10cm 6cm},clip]{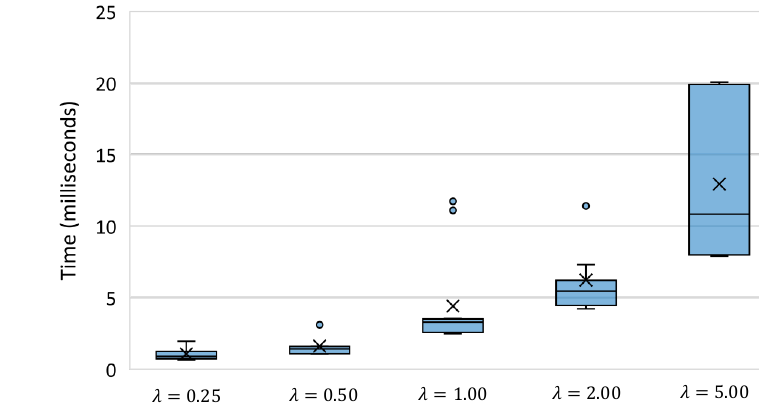}
    }
    \subfigure[$k=3$]{
    \includegraphics[width=0.45\textwidth,trim={11cm 6cm 10cm 6cm},clip]{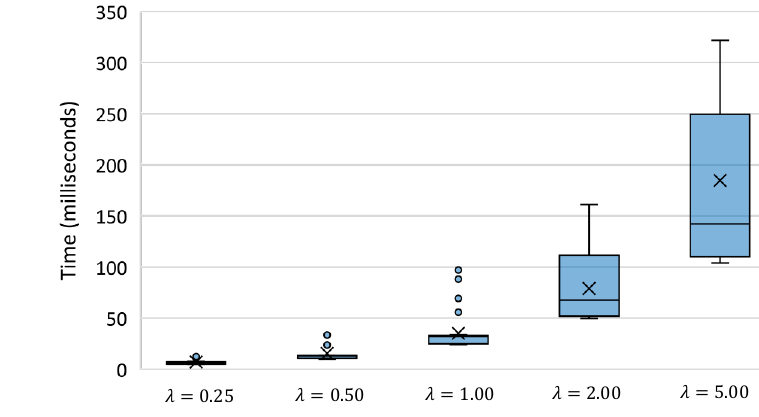}
    }
    
    \caption{Time to solve monitoring problems \eqref{eq:inference} with $n=200$ for different bandwidths $k$ and smoothness parameters $\lambda$, in milliseconds. Each plot depicts distributions over different signals and sparsity parameters $\mu$. Setup costs used to compute the convex hulls, which can be done offline a priori, are under one minute for $k=2$ and under 15 minutes for $k=3$.}
    \label{fig:timesOnline}
    \end{center}
\end{figure}

 The second idea relates to the methodology used to compute the convex hulls. Time-series problems are typically amenable to dynamic programming algorithms, and we rely on similar ideas to convexify the sets. Throughout the paper, we express the methods using \emph{decision diagrams}. As a by-product of the chosen methodology, we find that problem \eqref{eq:inferenceConvex} can be directly solved as a shortest path problem on a directed acyclic graph. This approach allow us to further reduce the computational times from a few seconds to milliseconds. Moreover, we show that the size of the diagrams can be bounded by a quantity that is polynomial in $n$ and a precision parameter $\varepsilon$, resulting in a fully polynomial time approximation scheme (FPTAS) for problem \eqref{eq:miqo} when the bandwidth and condition number of $\bm{Q}$ are fixed and there are no side constraints \eqref{eq:miqo_constraints}. Figure~\ref{fig:timesOnline} summarizes some runtimes from experiments in our computational section, where problem \eqref{eq:inference} is solved to optimality with the moving average filter \eqref{eq:movingAverage} and $n=200$ using financial time-series data. In the experiments, the convex hull of $X_Q$ is computed once, a process requiring 15 minutes at most and often much less, and then used to sequentially solve close to 7,000 instances of \eqref{eq:inference} to optimality; the time reported corresponds to the time per instance.

The proposed approach has three additional ``hidden" benefits. First, while using off-the-shelf solvers is time and memory intensive, shortest path problems in acyclic graphs are easier to implement and deploy. In particular, once the computation of the convex hull is done, inference problem \eqref{eq:inference} can then be solved in edge devices in sensor networks, wearables, or other devices with limited CPU and RAM. Second, similar to the data $\bm{y}$, the parameter $\mu$ in \eqref{eq:inference} does not affect set $X_Q$: as a consequence, cross-validation to select this parameter can be performed much more efficiently, even in the offline setting. Third, the proposed method can accommodate additional constraints \eqref{eq:miqo_constraints}, provided that they can be conveniently represented using decision diagrams, allowing for the inclusion of additional priors in inference problems \eqref{eq:inference}.

 \subsection{Outline and notation}

 The rest of this paper is organized as follows. In \S\ref{sec:background} we present relevant background on decision diagrams, which have been used for the most part for MILO problems. In \S\ref{sec:dd} we present the main method which uses decision diagrams to solve \eqref{eq:inference} and, more generally, problem \eqref{eq:miqo} with banded matrices. In \S\ref{sec:bandedMatrix} we show that the proposed method can be modified into an FPTAS for problem \eqref{eq:miqo_objective}-\eqref{eq:miqo_indicators} over matrices with fixed bandwidth and condition number. In \S\ref{sec:convexification} we formally establish the connection between the decision diagrams produced by the proposed method and the convex hull of set $X_Q$. Finally, \S\ref{sec:computations} presents computational experiments.

\subsubsection*{Notation} Throughout the paper, we denote vectors and matrices in \textbf{bold}. Given $n\in \Z_+$, let $[n]\defeq\left\{1,\dots,n\right\}$. Given a vector $\bm{v}\in \R^n$ or matrix $\bm{V}\in \R^{n\times n}$ and a set $S\subseteq [n]$, we denote by $\bm{v_S}$ and $\bm{V_S}$ the subvector or principle submatrix induced by $S$, respectively. Given two matrices $\bm{V}$ and $\bm{W}$ of same dimensions, we denote by $\langle \bm{V},\bm{W}\rangle\defeq \sum_{i}\sum_j V_{ij}W_{ij}$ the inner product between the matrices, and by $\bm{V}\circ \bm{W}$ their Hadamard (entrywise) product, that is, $(\bm{V}\circ \bm{W})_{ij}=V_{ij}W_{ij}$. Given $n,m\in \Z_+$, we denote by $\bm{0}_n$ an $n$-dimensional vector of zeros and by $\bm{0}_{n\times m}$ an $n\times m$-dimensional matrix of zeros; when the dimensions are clear from the context, we omit the subscripts. Given $S\subseteq [n]$, we let $\bm{e_S}\in \{0,1\}^n$ denote the vector that has ones in the positions indexed by $S$ and $0$ elsewhere. Moreover, given $i\in [n]$, we let $\bm{e_i}\defeq\bm{e_{\{i\}}}$ denote the $i$-th basis vector of $\R^n$. Finally, given a matrix $\bm{W}\in \R^{n\times n}$, we let $\bm{W_j}\in \R^n$ denote the $j$-th column of $\bm{W}$. 

This paper extensively uses pseudoinverses. 
Given a matrix $\bm{W}$, we let $\bm{W}^\dagger$ denote its Moore-Penrose pseudoinverse.  A special case of the pseudoinverse that is used throughout the paper pertains to matrices of the form $\bm{W}=\begin{pmatrix}\bm{A}&\bm{0}\\
\bm{0}&\bm{0}\end{pmatrix}$ where $\bm{A}$ is invertible, in which case $\bm{W}^\dagger=\begin{pmatrix}\bm{A^{-1}}&\bm{0}\\
\bm{0}&\bm{0}\end{pmatrix}$.  Given a matrix $\bm{W}\in \R^{n\times n}$ and set $S$, denote by $\bm{W_S}\in \R^{S\times S}$ the submatrix of $\bm{W}$ induced by $S$. Moreover, given $\bm{W}\in \R^{n\times n}$ and $S\subseteq [n]$,  observe that $\left(\bm{W}\circ \bm{e_Se_S^\top}\right)\in \R^{n\times n}$ is the matrix that coincides with $\bm{W_S}$ in the entries indexed by $S$ and is $0$ elsewhere. Similarly, if $\bm{W_S}$ is invertible, then $\left(\bm{W}\circ \bm{e_Se_S^\top}\right)^\dagger\in \R^{n\times n}$ is the matrix that coincides with $\bm{W_S^{-1}}$ in the entries indexed by $S$ and is $0$ elsewhere.

\begin{example}
Consider matrix $\bm{W}=\begin{pmatrix}2&1\\ 1&3\end{pmatrix}$. Then we find that 
\footnotesize\begin{align*}
\left(\bm{W}\circ \bm{0}\bm{0}^\top\right)^\dagger =\begin{pmatrix}0&0\\
0&0\end{pmatrix},\; \left(\bm{W}\circ \bm{e_1}\bm{e_1}^\top\right)^\dagger =\begin{pmatrix}1/2&0\\
0&0\end{pmatrix},\; \left(\bm{W}\circ \bm{e_2}\bm{e_2}^\top\right)^\dagger =\begin{pmatrix}0&0\\
0&1/3\end{pmatrix},\; \left(\bm{W}\circ \bm{e_{\{1,2\}}}\bm{e_{\{1,2\}}}^\top\right)^\dagger =\begin{pmatrix}3/5&-1/5\\
-1/5&2/5\end{pmatrix}\cdot\hfill \blacksquare
\end{align*}\normalsize
\end{example}

\section{Background on decision diagrams}\label{sec:background}

Decision diagrams are fundamental to the algorithm we propose in the paper. Our work departs from the classical literature in decision diagrams in the sense that we use them to tackle a nonlinear mixed-integer optimization problem and construct a non-polyhedral convex hull; \S\ref{sec:dd} is devoted to explaining our algorithm. Before, we review fundamental concepts concerning decision diagrams in linear settings. These concepts are necessary to model constraints \eqref{eq:miqo_constraints} when $Z\neq \{0,1\}^n$.

\subsection{Decision Diagrams for MILOs: Definitions and Notation}\label{sec:ddInMilo}
 A binary decision diagram $\dd = (\Nodes,\A,\valvec,\lenvec)$ encodes a set of points $\feasibleSet \subseteq \{0,1\}^n$ as a directed acyclic graph with node set $\Nodes$ and arc set $\A \subseteq \Nodes \times \Nodes$. Set $\Nodes$ is partitioned into $n+1$ layers $\Nodes_1, \dots, \Nodes_{n+1}$, where $\Nodes_1 = \{\rootnode\}$ contains a \emph{root} node $\rootnode$ and $\Nodes_{n+1}$ contains a set of \emph{terminal nodes}. Arcs $a \in \A$ connect nodes in consecutive layers and represent value assignments for the components of $\bm{z} \in \feasibleSet$. Let $\ell(a) \in \{1,\dots,n\}$ be the layer from which arc $a \in \A$ emanates. We denote the tail node of an arc by $\tail{a} \in \Nodes_{\ell(a)}$, the head node by $\head{a} \in \Nodes_{\ell(a)+1}$, and its value assignment by $\val{a} \in \{0,1\}$. Typically, in the context of MILOs, each arc is also assigned a length $l_a\in \R$ representing linear costs associated with assignment $\nu_a$.        

A decision diagram $\dd$ encodes $\feasibleSet$ through $\rootnode-\terminalnode$ paths, with $t \in \Nodes_{n+1}$, as follows. Each component of $\bm{z} \in \feasibleSet$ corresponds to a layer in $\dd$. An arc-specified $\rootnode-\terminalnode$ path $(a_1, \dots, a_n)$, where $\head{a_i} = \tail{a_{i+1}}$ for $i=1,\dots,n-1$, encodes the vector $\bm{z} = (\val{a_1}, \dots, \val{a_n})^\top$. A decision diagram is \emph{exact} if every point $\bm{z} \in \feasibleSet$ maps to a corresponding $\rootnode-\terminalnode$ path and vice versa. In the context of minimizing an objective function $f(\bm{z})$, an exact decision diagram satisfies that for any $\rootnode-\terminalnode$ path $(a_1, \dots, a_n)$ the path length $\sum_{j=1}^n \len{a_j} = f(\bm{z})$. As a result, a shortest path in an exact decision diagram $\dd$ yields an optimal solution to $\min_{\bm{z} \in \feasibleSet} \left\{ f(\bm{z}) \right\}$.

Decision diagrams are commonly generated from the state-transition graph of a dynamic programming (recursive) formulation of the problem \citep{bergman2016}. Such formulations consist of a \emph{state space}, a set of \emph{transition functions}, and a set of \emph{cost functions}. A dynamic program sequentially sets values for the decision variables, storing the outcome of these decisions in states, which store the relevant information about the partial solution obtained after fixing a subset of the variables. Transition functions establish how the system transitions between states, incurring a cost given by the cost functions. 

To illustrate the generation of a decision diagram from a state-transition graph, we consider the time series problem described above and a constraint that requires non-zero values to come in batches of at least $\tau$ consecutive periods. We refer to these constraints as \emph{contiguity constraints} throughout the paper. Letting $\bm{\zeta}\in \{0,1\}^n$ be auxiliary binary variables used to indicate whether a batch of consecutive non-zeros starts at time period $i$, set $Z$ is the set of binary points satisfying the linear inequalities
\begin{subequations}\label{eq:consecutiveOnes}
\begin{align}
&z_1\leq \zeta_1 \label{eq:consecutiveFirst} \\ 
&z_i-z_{i-1}\leq \zeta_i &\forall i\in \{2,\dots,n+1-\tau\}\\
&z_i-z_{i-1}\leq 0 &\forall i\in \{n+2-\tau,\dots,n\}\\
&\zeta_i\leq z_{i+j-1}&\forall  i\in \{1,\dots,n+1-\tau\},\; \forall j\in \{1,\dots,\tau\}. \label{eq:consecutiveLast}
\end{align}
\end{subequations}
To encode contiguity constraints using decision diagrams, we define the state variable $s^\ell$ as the minimum between $\tau$ and the number of consecutive non-zero values at decision stage $\ell$. Note that $s^\ell$ records the state of the system after having assigned values to $\ell-1$ variables. The initial state is given by $\{0\}$ and we use $\{\}$ to represent the infeasible state. We define the transition function as $s^{\ell+1} = \phi(s^{\ell},\hat{z}_{\ell})$, where $\hat{z}_{\ell}$ is the value assigned to variable $z_{\ell}$, and           
\begin{align} \label{eq:transitionsExample}
\phi(s^{\ell},\hat{z}_{\ell})=&
\begin{cases}
\{\} \quad \text{if } \hat{z}_{\ell} = 0 \text{ and } 1\leq s^\ell < \tau  \\
\{\} \quad \text{if } \ell = n, \; \hat{z}_{\ell} = 1, \text{ and } s^\ell < \tau-1  \\
\{0\} \quad \text{if } \hat{z}_{\ell} = 0 \text{ and } s^\ell = \tau  \\
\left\{ \min\{\tau,s^\ell+\hat{z}_{\ell}\} \right\} \quad \text{otherwise}
\end{cases}
\end{align}

The first case of \eqref{eq:transitionsExample} corresponds to the event in which at stage $\ell$, the number of consecutive non-zero values is greater than or equal to $1$ but less than $\tau$. In this case, it is not possible to set $z_\ell = 0$. The second case considers the event in which it is not possible to complete a batch of $\tau$ consecutive non-zeros at the last decision stage. The third case considers the event in which $z_\ell$ is set to zero, after a batch of at least $\tau$ consecutive non-zero values. In this case, the count of consecutive non-zero values resets back to $0$. The last case considers all other events. 

\begin{figure}[!h]
    \centering
    \subfigure[]{
    \includegraphics[width=0.3\textwidth,trim={3cm 20.4cm 11cm 2.5cm},clip]{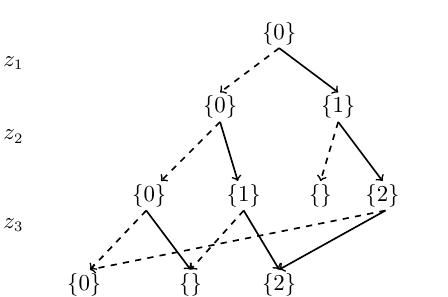}
    }
    \subfigure[]{
    \includegraphics[width=0.15\textwidth,trim={9cm 14.8cm 9cm 8.8cm},clip]{./figures/FigDD.pdf}
    }
\caption{State transition graph and corresponding decision diagram for a constraint that requires consecutive non-zero values.} 
    \label{fig:DP} 
\end{figure}

Figure \ref{fig:DP}(a) presents the state transition graph for the dynamic program described above with $n=3$ and $\tau = 2$. The corresponding decision diagram in Figure \ref{fig:DP}(b) is obtained in this example by removing the infeasible states. Note that there is a one-to-one correspondence between paths in the diagram and vectors that satisfy constraints~\eqref{eq:consecutiveOnes}. 
Solution approaches based on decision diagrams have shown promising results in discrete problems for which continuous relaxations tend to be weak  \citep{CirHoe13,BerCirHoeHoo14,BerCirHoeHoo16}). Recently, decision diagrams have also been used as convexification devices to obtain a characterization of the convex hull of a discrete feasible set in a lifted space. This approach is explored by \citet{lozano2018Binary} and \citet{macneil2024leveraging} to generate a convex representation of a discrete recourse problem in the context of two-stage stochastic programming. In the context of nonlinear optimization, \citet{davarnia2020} use decision diagrams to generate cuts for an outer approximation algorithm. We refer the reader to \citet{castro2022decision} for a comprehensive literature review on decision diagrams for optimization.

\subsection{Towards building decision diagrams for MIQO}

Often, in binary problems, the state space of the underlying dynamic program is conceptually simple and has a clear interpretation in terms of the problem to be solved. In the proposed approach, which involves both continuous variables and nonlinear terms involving the continuous variables, we use a state space that is perhaps less interpretable: selected columns of pseudoinverses associated with submatrices of $\bm{Q}$. We now provide a high-level intuition into why our selection of the state space is pertinent when solving \eqref{eq:miqo}, and defer to \S\ref{sec:dd} the formal statements.

Consider problem \eqref{eq:miqo} after projecting out the auxiliary variable $x_0$,
\begin{equation}\label{eq:optUnconstrained}\min_{\substack{\bm{x}\in \R^n,\bm{z}\in Z\\\bm{x}\circ(\bm{1}-\bm{z})=\bm{0}}}\bm{d^\top x}+\bm{c^\top z}+\frac{1}{2}\bm{x^\top Q x}.\end{equation}
We now proceed to write \eqref{eq:optUnconstrained} as a purely combinatorial optimization problem by projecting out the continuous variables $\bm{x}$. Indeed, given any $\bm{z}\in \{0,1\}^n$, and letting $S=\left\{i\in [n]:z_i=1\right\}$, it follows from the optimality conditions for the continuous variables that $\bm{x_S^*}=-\bm{Q_S^{-1}d_S}$ and $\bm{x_{[n]\setminus S}^*}=\bm{0}$, where $\bm{x_S^*}$ and $\bm{x_{[n]\setminus S}^*}$ denote the subvectors of $\bm{x^*}$ induced by entries in $S$ and not in $S$, respectively. Thus, projecting out variable $\bm{x}$, we find that \eqref{eq:optUnconstrained} is equivalent to 
\begin{equation}\label{eq:optCombinatorial}\min_{S}-\frac{1}{2}\bm{d_S^\top Q_S^{-1}d_S}+\sum_{i\in S}c_i=\min_{\bm{z}\in Z}-\frac{1}{2}\bm{d^\top} \left(\bm{Q}\circ \bm{zz^\top}\right)^{\dagger}\bm{d}+\bm{c^\top z} 
\end{equation}
The equality in \eqref{eq:optCombinatorial} holds since entries of $\bm{d}$ associated with variables such that $z_i=0$ are multiplied by zero entries of matrix $\left(\bm{Q}\circ \bm{zz^\top}\right)^{\dagger}$. Since the objective term in \eqref{eq:optCombinatorial} is a linear function of matrix $\left(\bm{Q}\circ \bm{zz^\top}\right)^{\dagger}$, it suggests that understanding the properties of this matrix is critical to the solution of \eqref{eq:optUnconstrained}. We formalize this intuition in the next section.

\section{Decision diagrams for MIQO}\label{sec:dd}

In this section we discuss how to construct decision diagrams in order to solve MIQO problems. Initially, in \S\ref{sec:unconstrained}, we discuss a ``natural" construction of a decision diagram that, although general, is impractical as it always produces a diagram with exponentially many nodes and arcs. Then, in \S\ref{sec:sparseMatrix} we propose a refinement of the state space, which allows for a reduction in the size of the diagrams when matrix $\bm{Q}$ is banded. Finally, in \S\ref{sec:epsExact} we discuss considerations associated with numerical precision, which pave the way to the FPTAS we study in \S\ref{sec:bandedMatrix}.

As stated above, our choice of state space corresponds to pseudoinverses related to matrix $\bm{Q}$. The transition function of the dynamic program  thus needs to compute new pseudoinverse based on a existing state. Lemma~\ref{lem:blockwiseInversion} below and, more importantly, the ensuing corollary, are critical to the definition of the transition functions.

\begin{lemma}[Blockwise inversion, \citet{lu2002inverses}]~\label{lem:blockwiseInversion}
	Given a non-singular symmetric square matrix $\bm V=\begin{pmatrix}
    \bm A&\bm B^\top\\
    \bm B&\bm G
\end{pmatrix}$, its inverse is given by $$\bm V^{-1}
=\begin{pmatrix}
	\bm A^{-1}+ \bm A^{-1}\bm B^\top(\bm V/\bm A)^{-1}\bm B\bm A^{-1}&-\bm A^{-1}\bm B^\top(\bm V/\bm A)^{-1}\\ -(\bm V/\bm A)^{-1}\bm B\bm A^{-1} & (\bm V/\bm A)^{-1}
\end{pmatrix},
$$
 where $\bm V/\bm A\defeq\bm G-\bm B \bm{A^{-1}}\bm B^\top$ is the Shur complement of the block $\bm G$ in $\bm V$.
\end{lemma}

\begin{corollary}\label{cor:inversion}
	Given a positive definite matrix $\bm{V}= \begin{pmatrix}\bm{A} & \bm{v}\\\bm{v^\top}&\delta\end{pmatrix}$, where $\bm{A}\in \R^{n\times n}$, $\bm{v}\in \R^n$ and $\delta \in \R_{+}$, its inverse is given by \small$$\bm{V^{-1}}= \begin{pmatrix}\bm{A^{-1}}+\bm{A^{-1}}\bm{v}(\bm V/\bm A)^{-1}\bm{v^\top A^{-1}} & -\bm{A^{-1}v}(\bm V/\bm A)^{-1}\\(\bm V/\bm A)^{-1}\bm{v^\top A^{-1}}&(\bm V/\bm A)^{-1}\end{pmatrix}.$$\normalsize
	Moreover, letting $\bm{u}\defeq\begin{pmatrix}-\bm{A^{-1}v}\\1\end{pmatrix}$, the identity 
	\begin{equation}\label{eq:diffInverses}\begin{pmatrix}\bm{A} & \bm{v}\\\bm{v^\top}&\delta\end{pmatrix}^\dagger =\begin{pmatrix}\bm{A} & \bm{0}\\\bm{0^\top}&0\end{pmatrix}^\dagger+\frac{1}{\delta-\bm{v^\top A^{-1}v}}\bm{uu^\top}\end{equation}
	holds.
\end{corollary}

\subsection{Definition of full decision diagrams for unconstrained MIQOs}\label{sec:unconstrained}

For simplicity, and to focus on the task of modeling the nonlinear terms in \eqref{eq:optUnconstrained}, we initially assume that $Z=\{0,1\}^n$ and thus the MIQO reduces to
\begin{equation}\label{eq:optRepeated}\min_{\substack{\bm{x}\in \R^n,\bm{z}\in \{0,1\}^n\\\bm{x}\circ(\bm{1}-\bm{z})=\bm{0}}}\bm{d^\top x}+\bm{c^\top z}+\frac{1}{2}\bm{x^\top Q x}=\min_{\bm{z}\in \{0,1\}^n}-\frac{1}{2}\bm{d^\top} \left(\bm{Q}\circ \bm{zz^\top}\right)^{\dagger}\bm{d}+\bm{c^\top z}.\end{equation}
Nonetheless, as Remark~\ref{rem:constraints} in this section shows, other feasible sets $Z$ can be handled by the decision diagrams as well. 
 We now discuss a state variable and transition function to solve \eqref{eq:optRepeated} to optimality, where the variable ordering is given by the natural order stemming from the underlying time series.  

We define state variables to be $n\times n$ matrices, representing matrices of the form $\left(\bm{Q}\circ \bm{ \bar z \bar z^\top}\right)^\dagger$ where $\bm{\bar z}\in \{0,1\}^n$ is a vector denoting the partial assignment of the indicator variables. More formally, we define the initial state $s^1=\left\{\bm{0}_{n\times n}\right\}$, and we define the transition function as follows: given state $s^\ell=\left\{\bm{\bar W^\ell}\right\}$, let vector $\bm{u}\in \R^n$ be \begin{equation}\label{eq:defU}\bm{u}=\frac{1}{\sqrt{Q_{\ell\ell}-\bm{Q_\ell^\top \bar W^\ell Q_\ell}}}\left(-\bm{\bar W^\ell}\bm{Q_\ell}+ \bm{e_\ell}\right)=\frac{1}{\sqrt{Q_{\ell\ell}-\sum_{i=1}^n\sum_{j=1}^n\bar{W}_{ij}^\ell Q_{i\ell}Q_{j\ell}}}\left(-\sum_{i=1}^n\bm{\bar W_i^\ell}Q_{i\ell}+ \bm{e_\ell}\right),\end{equation} let $s^{\ell+1} = \phi_{\text{full}}(s^{\ell},\hat{z}_{\ell})$, where $\hat{z}_{\ell}$ is the value assigned to variable $z_{\ell}$, and         
\begin{align} \label{eq:transitionsFull}
\phi_\text{full}(s^{\ell},\hat{z}_{\ell})=&
\begin{cases}
\left\{\bm{\bar W^\ell}\right\} \quad \text{if } \hat{z}_{\ell} = 0 \\
\left\{ \bm{\bar W^\ell}+\bm{uu^\top} \right\}\quad \text{if } \hat{z}_{\ell} = 1.
\end{cases}
\end{align}

\begin{definition}[Complete decision diagram]\label{def:fullDD}
Define the state transition graph $\mathcal{G}_{\text{full}}=(\Nodes,\A)$ by applying the state transition function $\phi_\text{full}$ recursively, starting from the initial state $s^1$. Define the complete decision diagram $\dd_{\text{full}}=(\mathcal{G}_{\text{full}},\valvec,\Uvec)$ as the resulting graph in which each node corresponds to a state, there is a directed arc $(s^i,\phi_{\text{full}}(s^i, \hat z))$ for  $\hat z\in \{0,1\}$, the value assignment function $\nu:A\to \{0,1\}$ is such that $\nu (s^i,\phi_{\text{full}}(s^i, \hat z))=\hat z$, and the transition vectors $\bm{u}:A\to \R^n$ are such that $\bm{u_a}=\bm{0}$ if $\nu(a)=0$ and are given by \eqref{eq:defU} if $\nu(a)=1$, where $\bm{\bar W^\ell}$ is the matrix stored by the state at the tail of the arc $\tail{a}$. 
\end{definition}

To simplify the notation, given any arc $a\in A$, we will use $\nu_a$ and $\bm{u_a}$ instead of $\nu(a)$ and $\bm{u}(a)$. We will also refer to $\nu_a$ and $\bm{u_a}$ as the value assignment and transition vector stored in arc $a$ of the diagram, respectively.

\begin{example} \label{ex:fullDD}Consider problem \eqref{eq:optRepeated} with matrix $\bm{Q}$ given by $\footnotesize \begin{pmatrix}2 & -1 & -1\\
-1 & 3 & -1\\
-1 & -1 & 2\end{pmatrix}.$ Figure~\ref{fig:Example1}(a) shows the state transition graph obtained by applying function $\phi_{\text{full}}$ recursively, and Figure~\ref{fig:Example1}(b) depicts the transition vectors stored in the arcs of the diagram.\hfill $\blacksquare$
\begin{figure}[!h]
    \centering
    \subfigure[]{
    \includegraphics[width=0.65\textwidth,trim={3cm 19cm 4.8cm 2.5cm},clip]{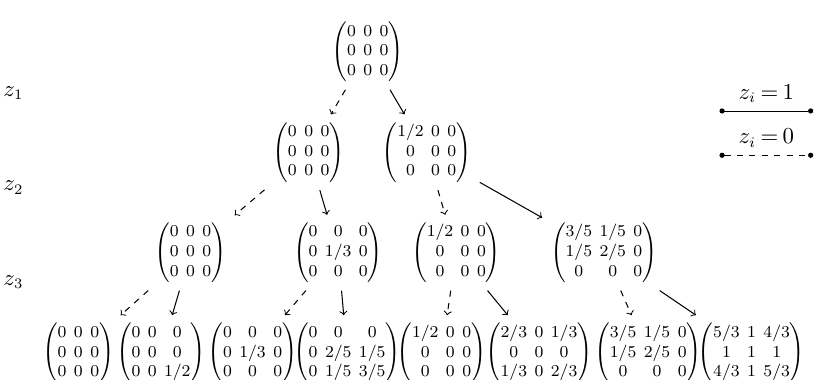}
    } \\
    \subfigure[]{
    \includegraphics[width=0.65\textwidth,trim={1cm 19.8cm 6.2cm 2.5cm},clip]{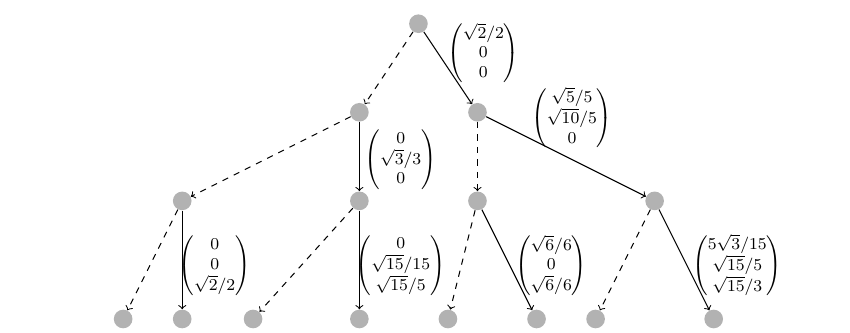}
    }
    
    \caption{State transition graph (a) and transition vectors (b) generated for Example~\ref{ex:fullDD}.} 
    \label{fig:Example1}
\end{figure}

\end{example}

Before formally proving the key properties of the decision diagram generated (Proposition~\ref{prop:state}), we state in Remarks \ref{rem:noLength}-\ref{rem:zeros} below some useful observations concerning states and transition vectors.

\begin{remark}\label{rem:noLength} Departing from standard practices, $\dd_{\text{full}}$ does not have associated arc lengths. Indeed, to allow for real-time execution, the diagram is constructed without observing the linear coefficients $\bm{c}$ and $\bm{d}$ in \eqref{eq:optRepeated}. Instead we store the transition vectors $\bm{u}$ which as we show later in Proposition~\ref{prop:pathLength}, are sufficient to compute lengths once the data $\bm{c},\bm{d}$ are realized. $\hfill \blacksquare$
\end{remark}

\begin{remark}\label{rem:uDiff}
From the definition of the state transition function \eqref{eq:transitionsFull}, given any arc $a\in A$ with tail and head states storing matrices given by $\bm{\bar W^\ell}$ and $\bm{\bar W^{\ell+1}}$, respectively, we find that $\bm{u_a u_a^\top}=\bm{\bar W^{\ell+1}}-\bm{\bar W^\ell}$. Moreover, vector $\bm{u}$ can also be retrieved directly from the $\ell$-th column of the head state as $\bm u_{a}=\bm{\bar W_{\ell}^{\ell+1}}/\sqrt{{\bar W}^{\ell+1}_{\ell\ell}}$
$\hfill \blacksquare$
\end{remark}

\begin{remark}\label{rem:zeros}
From Figure~\ref{fig:Example1} we observe that the sparsity patterns of the diagonal elements of states $\bm{W^\ell}$ is different for all nodes in the same layer. In fact, this sparsity pattern corresponds exactly to the sparsity pattern of the vector $\bm{\bar z}$ encoded by the path leading to each node. This property, which is not necessarily evident from the definition of the transition function, holds true in general (and is a direct corollary of Proposition~\ref{prop:state} below). As a consequence, since different paths in the diagram encode different vectors $\bm{\bar z}$, each with a unique sparsity pattern, it follows that the number of nodes in the $\ell$-th layer is always $2^{\ell-1}$ and $\dd_{\text{full}}$ corresponds to complete enumeration. Clearly, we do not advocate for using $\dd_{\text{full}}$ in practice.$\hfill \blacksquare$
\end{remark}

We now prove the fundamental property associated with $\dd_{\text{full}}$, namely, that the states generated by transition function $\phi_{\text{full}}$ indeed correspond to pseudo-inverses associated with matrix $\bm{Q}$.

\begin{proposition}\label{prop:state} Given any arc-specified path $(a_1,a_2,\ldots,a_{\ell-1})$ in $\mathcal{G}_{\text{full}}$ such that $\tail{a_1}=s^1$ and $\head{a_{\ell-1}}=s^\ell$, letting $\bm{\bar W^\ell}$ be the matrix stored by state $s^\ell$ and $\bm{\bar z^\ell}=\sum_{i=1}^{\ell-1}\nu_{a_i}\bm{e_i}$ be the partial solution represented by the path, the identity $\bm{\bar W^\ell}=\left(\bm{Q}\circ\bm{\bar z^\ell}(\bm{\bar z^\ell})^\top\right)^\dagger$ holds.
\end{proposition}
\proof
Observe that since $\bm{\bar z^\ell}$ corresponds to a partial solution with no assignments for variables $\{z_\ell, z_{\ell+1},\dots, z_n\}$, we have that $\bar z_i^\ell=0$ for $i\geq \ell$. We prove the result by induction on the layer $\ell$.

\noindent \textbf{Base case} If $\ell=1$, then $\bm{\bar z^\ell}=\bm{0}_n$ and $s^\ell=\bm{0}_{n\times n}$, and the result is trivially satisfied.

\noindent \textbf{Inductive step} Assume that the result holds for $s^\ell$. There are two cases, depending on the assignment value $\hat z_\ell$. \newline
$\bullet$ If $\hat z_\ell=0$, then the result is automatically satisfied since in this case the state does not change ($\bm{\bar W^{\ell+1}}=\bm{\bar W^\ell}$) and $\bm{\bar z^{\ell+1}}=\bm{\bar z^\ell}$.\newline
$\bullet$ If $\hat z_\ell=1$, then \begin{align*}\bm{u}&=\frac{1}{\sqrt{Q_{\ell\ell}-\sum_{i=1}^n\sum_{j=1}^n\bar{W}_{ij}^\ell Q_{i\ell}Q_{j\ell}}}\left(-\sum_{i=1}^n\bm{\bar W_i^\ell}Q_{i\ell}+ \bm{e_\ell}\right)\\
&=\frac{1}{\sqrt{Q_{\ell\ell}-\sum_{i=1}^n\sum_{j=1}^n\left(\bm{Q}\circ\bm{\bar z^\ell}(\bm{\bar z^\ell})^\top\right)_{ij}^\dagger Q_{i\ell}Q_{j\ell}}}\left(-\sum_{i=1}^n\left(\bm{Q}\circ\bm{\bar z^\ell}(\bm{\bar z^\ell})^\top\right)_i^\dagger Q_{i\ell}+ \bm{e_\ell}\right)\tag{$\because$ induction hypothesis}\\
&=\frac{1}{\sqrt{Q_{\ell\ell}-\sum_{i=1}^{\ell-1}\sum_{j=1}^{\ell-1}\left(\bm{Q}\circ\bm{\bar z^\ell}(\bm{\bar z^\ell})^\top\right)_{ij}^\dagger Q_{i\ell}Q_{j\ell}}}\left(-\sum_{i=1}^{\ell-1}\left(\bm{Q}\circ\bm{\bar z^\ell}(\bm{\bar z^\ell})^\top\right)_i^\dagger Q_{i\ell}+ \bm{e_\ell}\right)\tag{$\because \bar z_i^\ell=0$ for $i\geq \ell\implies \left(\bm{Q}\circ\bm{\bar z^\ell}(\bm{\bar z^\ell})^\top\right)_{ij}^\dagger=0$ if $\max\{i,j\}\geq \ell$}.
\end{align*}
In particular, from Corollary~\ref{cor:inversion},we find that $\bm{uu^\top}=\left(\bm{Q}\circ\bm{\bar z^{\ell+1}}(\bm{\bar z^{\ell+1}})^\top\right)^\dagger-\left(\bm{Q}\circ\bm{\bar z^\ell}(\bm{\bar z^\ell})^\top\right)^\dagger$, where $\bm{\bar z^{\ell+1}}=\bm{\bar z^\ell}+\bm{e_\ell}$. The identity
\begin{align*}
\bm{\bar W^{\ell+1}}=\bm{\bar W^\ell}+\bm{uu^\top}&=\left(\bm{Q}\circ\bm{\bar z^\ell}(\bm{\bar z^\ell})^\top\right)^\dagger +\bm{uu^\top}=\left(\bm{Q}\circ\bm{\bar z^{\ell+1}}(\bm{\bar z^{\ell+1}})^\top\right)^\dagger
\end{align*}
follows immediately. 
\endproof

Finally, Proposition~\ref{prop:pathLength} shows how to solve problem \eqref{eq:optUnconstrained} as a shortest path in the (exponentially large) decision diagram $\dd_{\text{full}}$ once the data $\bm{c},\bm{d}$ are realized.
\begin{proposition}\label{prop:pathLength} Given vectors $\bm{c},\bm{d}\in \R^n$ and a full decision diagram $\dd_\text{full}$, define the length of an arc $a\in A$ as $l_a=c_{\tail{a}}\nu_{\tail{a}}-\frac{1}{2}\left(\bm{d^\top u_{a}}\right)^2.$ Then the length of any arc-specified path $(a_1,a_2,\ldots,a_{\ell-1})$ in $\mathcal{G}_{\text{full}}$  is given by 
	$$h(\bm{\bar z^\ell})\defeq\bm{c^\top \bar z^\ell}+ \min_{\substack{\bm{x}\in \R^n\\\bm{x}\circ(\bm{1}-\bm{\bar z^\ell})=\bm{0}}}\left\{\bm{d^\top x}+\frac{1}{2}\bm{x^\top Q x}\right\}=\bm{c^\top \bar z^\ell}-\frac{1}{2}\bm{d^\top} \left(\bm{Q}\circ \bm{\bar z^\ell}(\bm{\bar z^\ell})^\top\right)^{\dagger}\bm{d},$$
 where $\bm{\bar z^\ell}=\sum_{i=1}^{\ell-1}\nu_{a_i}\bm{e_i}$ is the partial solution represented by the path. In particular, any shortest path between the root and a terminal node corresponds to an optimal solution of~\eqref{eq:optUnconstrained}.
\end{proposition}
\proof
We do the proof by induction on the layer $\ell$.

\noindent \textbf{Base case $\bm{\ell=1}$} In this case $\bm{\bar z^1}=\bm{0}_n$. Defining the empty path as having length $0$, the result holds automatically. 

\noindent \textbf{Inductive step} We assume the result is true for path $(a_1,a_2,\dots,a_{\ell-2})$, corresponding to partial solution $\bm{\bar z^{\ell-1}}$ and where state $\head{a_{\ell-2}}$ stores matrix $\bm{\bar W^{\ell-1}}$, and prove it for the path including arc $a_{\ell-1}$, depending on the assignment $\nu_{a_{\ell-1}}$.

\noindent \textbf{Case $\bm{\nu_{a_{\ell-1}}=0}$} In this case $\bm{u}_{a_{\ell-1}}=\bm{0}_n$, and it follows that $l_{a_{\ell-1}}=0$. Moreover, since $\bm{\bar z^\ell}=\bm{\bar z^{\ell-1}}$, the result holds by the inductive hypothesis.  

\noindent \textbf{Case $\bm{\nu_{a_{\ell-1}}=1}$} In this case $\bm{\bar z^{\ell}}=\bm{\bar z^{\ell-1}}+\bm{e_{\ell-1}}$. In order to simplify the notation, we denote by $(\bm{Q^\ell})^\dagger =\left(\bm{Q}\circ \bm{\bar z^\ell}(\bm{\bar z^\ell})^\top\right)^\dagger.$ We find that
\begin{align*}h(\bm{\bar z^{\ell}})&=\bm{c^\top \bar z^{\ell-1}}+c_{\ell-1}-\frac{1}{2}\bm{d^\top}(\bm{Q^\ell})^\dagger\bm{d}\\
	&=\bm{c^\top \bar z^{\ell-1}}-\frac{1}{2}\bm{d^\top} (\bm{Q^{\ell-1}})^\dagger\bm{d}+c_{\ell-1}-\frac{1}{2}\bm{d^\top} \left((\bm{Q^\ell})^\dagger-(\bm{Q^{\ell-1}})^\dagger\right)\bm{d}\\
	&=h(\bm{\bar z^{\ell-1}})+c_{\ell-1}-\frac{1}{2}\bm{d^\top} \left((\bm{Q^\ell})^\dagger-(\bm{Q^{\ell-1}})^\dagger\right)\bm{d}\tag{$\because$ Induction hypothesis}\\
 &=h(\bm{\bar z^{\ell-1}})+c_{\ell-1}-\frac{1}{2}\bm{d^\top} \left(\bm{u_{a_{\ell-1}}u_{a_{\ell-1}}^\top}\right)\bm{d}\tag{$\because$ Remark~\ref{rem:uDiff}}\\
 &=h(\bm{\bar z^{\ell-1}})+c_{\ell-1}-\frac{1}{2}\left(\bm{d^\top u_{a_{\ell-1}}}\right)^2=h(\bm{\bar z^{\ell-1}})+l_{a_{\ell-1}},
\end{align*}
concluding the proof. 
\endproof

\subsection{Compressing the state space}\label{sec:sparseMatrix}

The dynamic program discussed in \S\ref{sec:unconstrained} is valid for any $\bm{Q}\succ \bm{0}$, but does not exploit any properties of the matrix such as the bandwidth, resulting in an exponential growth of the size of the state transition graph. 
We now discuss how to refine the state of the dynamic program for settings where matrix $\bm{Q}$ is a banded matrix, ultimately leading to a practical algorithm.

Recall the transition function defined in \eqref{eq:defU}-\eqref{eq:transitionsFull}. Observe that, that if $Q_{i\ell}=0$, then the $i$-th column of $\bm{\bar W^\ell}$ is \emph{irrelevant} to the computation of $\bm{u}$ in \eqref{eq:defU}. Moreover, if $\bm{Q}$ has bandwidth $k$, then $Q_{i\ell}=0$ whenever $\ell-i>k$, in other words, column $i$ of the state variable can be safely removed from the state for layers $\ell>k+i$. We now present the notion of relevant index, corresponding to the last layer in which a given column is relevant.

\begin{definition}[Relevant index]
Given matrix $\bm{Q}$ and index $i\in [n]$, define the \emph{relevant index} $\pi_i\in [n]$ as $$\pi_i\defeq\max_{j\in [n]}\{j: Q_{ij}\neq 0\}.$$ 
\end{definition}
Since we assume that $\bm{Q}\succ 0$ and therefore $Q_{ii}>0$ for all $i\in [n]$, we find that $i\leq \pi_i$. Moreover, if $\bm{Q}$ is a banded matrix with bandwidth $k$, then $\pi_i\leq i+k$ for all $i\in [n]$. Since $Q_{i\ell}=0$ for all $\ell>\pi_i$, it follows that the $i$-column of $\bm{\bar W}^\ell$ is irrelevant to the computation of transition vectors $\bm{u}$ in \eqref{eq:defU} whenever $\ell>\pi_i$. Moreover, from Proposition~\ref{prop:state}, we can infer that columns of $\bm{\bar W^\ell}$ with indexes $i\geq \ell$ are zero columns, and hence are irrelevant as well. Thus, given a state $\bm{\bar W}^\ell$, we say that the \emph{relevant columns} of $\bm{\bar W}^\ell$ are precisely the columns corresponding to indexes such that $i<\ell\leq \pi_i$.

We now update the state and transition functions of the dynamic program. We let $s^1=\left\{\bm{0}_{n\times n}\right\}$ and, given state $s^\ell=\left\{\bm{\bar W^\ell}\right\}$ and assignment value $\hat z_\ell$, let $\bm{u}$ be defined as in \eqref{eq:defU} and let the new transition function be $s^{\ell+1} = \phi(s^{\ell},\hat{z}_{i}) = \left\{ \bm{\bar W^{\ell+1}}\right\}$, where 
\begin{align} \label{eq:transitionsSparse}
\bar W_{ij}^{\ell+1}=&
\begin{cases}
0&\text{if }\ell\geq \pi_j\\
\bar W_{ij}^\ell & \text{if }\ell< \pi_j\text{ and } \hat{z}_{\ell} = 0 \\
\bar W_{ij}^\ell+u_iu_j & \text{if }\ell< \pi_j\text{ and } \hat{z}_{\ell} = 1.
\end{cases}
\end{align}
Observe that the first condition in \eqref{eq:transitionsSparse} sets previously computed columns of $\bm{\bar W^\ell}$ to zero whenever they are no longer relevant, while the other conditions match those in \eqref{eq:transitionsFull}.

\begin{definition}[Compressed Decision Diagram] \label{def:compressedDD} Define the state transition graph $\mathcal{G}=(\Nodes,\A)$ by applying the state transition function $\phi$ in \eqref{eq:transitionsSparse} recursively, starting from the initial state $s^1$. Define the compressed decision diagram $\dd=(\mathcal{G},\valvec,\Uvec)$ as the resulting graph in which each node corresponds to a state, there is a directed arc $(s^i,\phi(s^i, \hat z))$ for  $\hat z\in \{0,1\}$, the value assignment function $\nu:A\to \{0,1\}$ is such that $\nu (s^i,\phi(s^i, \hat z))=\hat z$, and the transition vectors $\bm{u}:A\to \R^n$ are such that $\bm{u_a}=\bm{0}$ if $\nu(a)=0$ and are given by \eqref{eq:defU} if $\nu(a)=1$, where $\bm{\bar W^\ell}$ is the matrix stored by~$\tail{a}$. 
\end{definition}

\begin{example}\label{ex:DD}
Consider problems with matrix $$\footnotesize \bm{Q}=\begin{pmatrix}4&-1&-1&0&0\\-1&4&0&-1&0\\-1&0&4&0&-1\\0&-1&0&4&-1\\
0&0&-1&-1&4\end{pmatrix},$$ which has bandwidth $k=2$. In this case the relevant indexes are given by $\bm{\pi^\top}=( 3\;\; 4\;\; 5\;\; 5\;\;5).$ Figure \ref{fig:Example3} shows the ensuing compressed decision diagram. Observe that, to depict the states, we simply show the relevant columns at each layer (corresponding to the $k$ most recently computed columns in this case), as the others are zero by definition. Moreover, we only show for each column the first $\ell-1$ rows, as remaining rows are zero as well. Observe that after $k$ consecutive zero assignments, the relevant columns revert to $\bm{0}_{n\times k}$ (a property that is easily shown to be true for any banded matrix), reducing the size of the decision diagram. Additional reductions occur for this decision diagram: for example, state $ \footnotesize \begin{pmatrix}0 & 0\\ 0&0\\
0.25&0\\
0 & 0\end{pmatrix}$ in layer $5$ is reached by two different paths, and both correspond to value assignments of one in layer $3$. In the end, we observe that the compressed decision diagram has 11 nodes in layer $\ell=5$, while the full decision diagram (Definition~\ref{def:fullDD}) would have 16. In the last layer, all terminal nodes all compressed into a single one automatically since no columns are relevant at this point. $\hfill \blacksquare$
\begin{figure}[!h]
    \centering
    \subfigure[]{
    \includegraphics[width=0.65\textwidth,trim={3cm 18.2cm 3.3cm 2.5cm},clip]{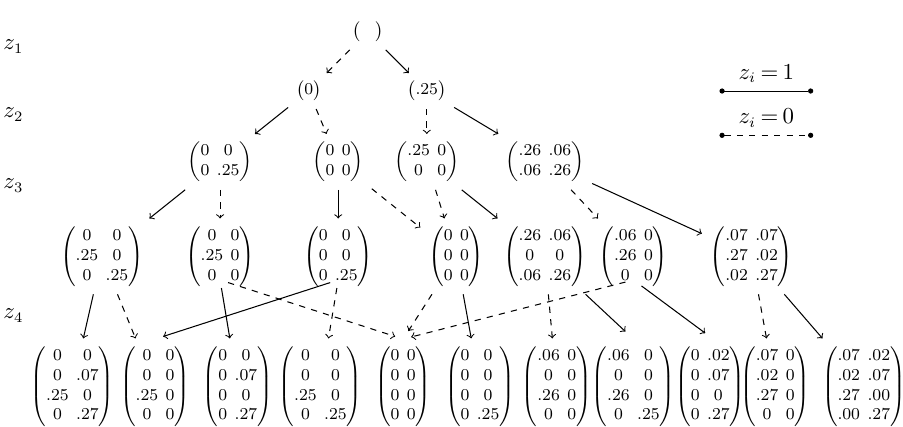}
    } \\
    \subfigure[]{
    \includegraphics[width=0.6\textwidth,trim={3cm 17.5cm 6.2cm 2.5cm},clip]{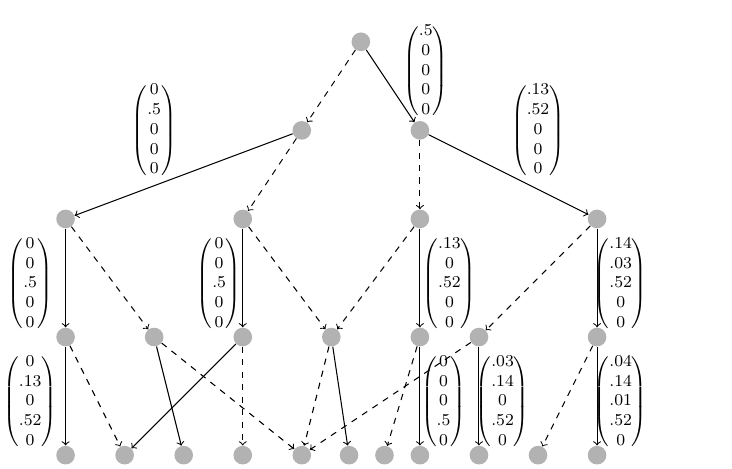}
    }
\caption{State transition graph (a) and transition vectors (b) generated for Example~\ref{ex:DD}. Last layer and some transition vectors omitted for clarity. In the graphs, ``0" denotes zero, while ``0.00" represents a number $0<\epsilon<0.005$. }\label{fig:Example3}
\end{figure}

\end{example}

\begin{proposition}\label{prop:correctReduced}
Given any arc-specified path $(a_1,a_2,\ldots,a_{\ell-1})$ in $\mathcal{G}$ such that $\tail{a_1}=s^1$ and $\head{a_{\ell-1}}=s^\ell$, letting $\bm{\bar W^\ell}$ be the matrix stored by state $s^\ell$ and $\bm{\bar z^\ell}=\sum_{i=1}^{\ell-1}\nu_{a_i}\bm{e_i}$ be the partial solution represented by the path,
the following properties hold true.
\begin{enumerate}
\item $\bar W_{ij}^\ell=\left(\bm{Q}\circ\bm{\bar z^\ell}(\bm{\bar z^\ell})^\top\right)^\dagger_{ij}$ if $\ell\leq\pi_j$ and $\bar W_{ij}^\ell=0$ if $\ell>\pi_j$
\item $\bm{u}_{a_{\ell-1}}=\left(\bm{Q}\circ\bm{\bar z^{\ell}}(\bm{\bar z^{\ell}})^\top\right)-\left(\bm{Q}\circ(\bm{\bar z^{\ell}}-\nu_{a_{\ell-1}}\bm{e_{\ell-1}})(\bm{\bar z^{\ell}}-\nu_{a_{\ell-1}}\bm{e_{\ell-1}})^\top\right)$.
\end{enumerate}
\end{proposition}

In the interest of completeness, we include a formal proof of Proposition~\ref{prop:correctReduced} in Appendix~\ref{sec:proofs}. A simple argument to show why the proposition holds is that compressed decision diagrams are obtained from the full diagrams by compressing nodes storing states that coincide in the relevant columns, and letting the matrix obtained by filling the irrelevant columns with zeros be the state stored in the compressed node. Since computations of vectors $\bm{u}$ are not affected by this procedure since either those columns were not used in the computations (if $\ell>\pi_i$) or the columns already contained zeros (if $\ell\leq i$), the rest of the decision diagram is not affected by this procedure. Therefore the compressed decision diagrams inherit the stated properties from the full diagrams. Similarly, since the vectors $\bm{u}$ stored in the arcs are identical to the ones in the full diagram, shortest path in compressed diagrams correspond to optimal solutions of \eqref{eq:optUnconstrained}. We state this fact succinctly in the next proposition.

\begin{proposition}\label{prop:pathLengthCompressed}
Proposition~\ref{prop:pathLength} holds true for compressed decision diagrams.
\end{proposition}

We close this section with three additional remarks concerning the implementation of the compressed decision diagrams, the special case arising from tridiagonal $\bm{Q}$ matrices, and the inclusion of additional constraints in the diagram.

\begin{remark}
While we chose to represent for mathematical convenience the state of compressed decision diagrams as a $n\times n$ matrix full of zero entries, as illustrated in Example~\ref{ex:DD}, a better implementation consists of storing only the relevant columns. In our code we use the more convenient approach of representing the matrix as $n\times n$ dimensional, but using libraries for sparse matrices to handle the mathematical operations. While this approach might result in a slight overhead versus only storing the relevant columns, it results in substantially better performance (especially memory-wise) than using standard dense matrices throughout the algorithm. $\hfill \blacksquare$
\end{remark}

\begin{remark}
In the special case of a tridiagonal matrix ($k=1$), where every arc $a$ with value assignments $\nu_a=0$ leads to the same $\bm{0}$ state, the number of nodes in each layer is $\ell+1$ and the total number of nodes produced is thus $\mathcal{O}(n^2)$. $\hfill \blacksquare$ 
\end{remark}

\begin{remark} \label{rem:constraints}
If constraints $\bm{z}\in Z$ can be easily represented via dynamic programming, they can be included in compressed decision diagrams as well by simply storing the state variable associated with the constraints. For example, to include contiguity constraints as discussed in \S\ref{sec:ddInMilo}, it suffices to include the number of consecutive nonzero values in the state, and compress states that coincide in the relevant columns of $\bm{\bar W^\ell}$ \emph{and} the value of consecutive nonzeros. $\hfill \blacksquare$
\end{remark}

\subsection{$\epsilon$-exact decision diagrams}\label{sec:epsExact}

Implementing the compressed decision diagrams discussed in \S\ref{sec:sparseMatrix} requires identifying whether two states in the same layer are equal. Recall that the states correspond to columns of inverses of submatrices of $\bm{Q}$. In theory, and assuming elements of $\bm{Q}$ are rational, it would be possible to perform all operations with exact arithmetic, thus accomplishing the compression exactly. Indeed, note that while the computation of $\bm{u}$ in \eqref{eq:defU} involves a square root, updates in \eqref{eq:transitionsFull} and \eqref{eq:transitionsSparse} require only matrix $\bm{uu^\top}$, hence the square root does not need to be computed. 
In practice, however, we use standard, finite-precision, libraries to compute sums, products and divisions, resulting in roundoff errors. As a consequence, an exact verification of whether two states are equal is no longer possible. Instead we turn to the concept of indistinguishable states, as defined next.

\begin{definition}[$\epsilon$-Indistinguishable states]\label{def:epsExact}
Given a decision diagram $\dd$, a parameter $\epsilon\geq 0$, and two states in the same layer $\bar s^\ell=\{\bm{\bar W^\ell}\}$ and $\hat s^\ell=\{\bm{\hat W^\ell}\}$, the states are $\epsilon$-indistinguishable if
$$\max_{j:j\leq \ell\leq \pi_j}\left\{\|\bm{\bar W_j^\ell}-\bm{\hat W_j^\ell}\|_\infty\right\} \leq \epsilon.$$
\end{definition}
Observe that the comparison in Definition~\ref{def:epsExact} is performed only for the relevant columns, as the remaining ones are zero in the compressed diagram. Interestingly, the definition can be used for complete decision diagrams $\dd_{\text{full}}$ as well. Moreover, two states are $0$-indistinguishable in $\dd_{\text{full}}$ if and only if the two solutions $\bar z^\ell$ and $\hat z^\ell$ represented by the unique paths leading to those states in $\mathcal{G}_{\text{full}}$ correspond to paths in the compressed decision diagram $\dd$ leading to a same state. In other words, a compressed diagram is obtained by ``merging" $0$-indistinguishable states in $\dd_{\text{full}}$. 
Therefore, merging $\epsilon$-indistinguishable states with $\epsilon>0$ would result in even smaller decision diagrams, at the cost of incurring approximation errors.

\begin{definition}[$\epsilon$-exact decision diagram]\label{def:ddEpsilon}
An $\epsilon$-exact decision diagram is any decision diagram produced layer by layer according to Definition~\ref{def:compressedDD} and \emph{merging} $\epsilon$-indistinguishable states (selecting any of the states in the merge as the representative). 
\end{definition}

Note that we do not explicitly distinguish whether the decision diagram is complete or compressed, as the ensuing $\epsilon$-exact diagram is the same (assuming comparisons between states are performed in the same order). While $\epsilon$-exact decision diagrams are approximations, we use the term ``exact" to emphasize that our original intention is not to artificially reduce their size, but simply to handle numerical imprecisions. In our computations in \S\ref{sec:computations}, we set $\epsilon=10^{-5}$, matching the default integrality tolerance of off-the-shelf  solvers, and recover in 99\% of the problems the same solution as off-the-shelf commercial solvers. In the remaining 1\% of the instances the maximum relative difference between the objective values reported by decision diagrams and mixed-integer solvers is at most 0.04\% -- and in some of these instances, the solution obtained from using decision diagrams is better. These results are presented in detail in Appendix~\ref{sec:numPrec}. In summary, $\epsilon$-exact decision diagrams can provide solutions as accurate as ``exact" MIQO solvers (which are also subject to numerical imprecisions) provided that $\epsilon$ is small enough. Nonetheless we also observed that for values of $\epsilon$ sufficiently small (e.g., $\epsilon=10^{-5}$ vs $\epsilon=10^{-6}$) the solutions obtained are essentially identical, but the size of the diagrams could vary substantially. 

\begin{example}\label{ex:merging}
Consider a problem with matrix {\footnotesize$ \bm{Q}=\begin{pmatrix}5 & -1 & 0 & 0 & 0&0&0\\
-1&5& -1& 0&0&0&0\\
0&-1&5&-1&0&0&0\\
0&0&-1&5&-1&0&0\\
0&0&0&-1&5&-1&0\\
0&0&0&0&-1&5&*\\
0&0&0&0&0&*&*\end{pmatrix},$} where $*$ are some unspecified values associated with variable $x_7$. Moreover, consider states at layer $\ell=7$, where only variable $z_7$ remains unassigned. 
Since the matrix is tridiagonal, $\pi_j=j+1$ and only the 6th column of the state variable is preserved in the compressed decision diagrams. Consider the partial solutions $(\bm{\bar z^\ell})^\top = (1\;\;1\;\;1\;\;1\;\;1\;\;1\;\;0)$ and $(\bm{\hat z^\ell})^\top = (0\;\;1\;\;1\;\;1\;\;1\;\;1\;\;0)$: we find that 
\footnotesize\begin{align*}\left(\bm{Q}\circ \bm{\bar z^\ell}(\bm{\bar z^\ell})^\top\right)^\dagger&=\begin{pmatrix}0.21&0.04&0.01&0.00&0.00&\textbf{0.00008}&0\\
0.04&0.22&0.05&0.01&0.00&\textbf{0.00040}&0\\
0.01&0.05&0.22&0.05&0.01&\textbf{0.00190}&0\\
0.00&0.01&0.05&0.22&0.05&\textbf{0.00909}&0\\
0.00&0.00&0.01&0.05&0.22&\textbf{0.04356}&0\\
0.00&0.00&0.00&0.01&0.04&\textbf{0.20871}&0\\
0&0&0&0&0&\textbf{0}&0\\
\end{pmatrix}\\
\left(\bm{Q}\circ \bm{\hat z^\ell}(\bm{\hat z^\ell})^\top\right)^\dagger&=\begin{pmatrix}0&0&0&0&0&\textbf{0}&0\\
0&0.21&0.04&0.01&0.00&\textbf{0.00038}&0\\
0&0.04&0.22&0.05&0.01&\textbf{0.00189}&0\\
0&0.01&0.05&0.22&0.05&\textbf{0.00909}&0\\
0&0.00&0.01&0.05&0.22&\textbf{0.04356}&0\\
0&0.00&0.00&0.01&0.04&\textbf{0.20871}&0\\
0&0&0&0&0&\textbf{0}&0\\
\end{pmatrix},\end{align*}\normalsize
where in both cases we use ``$0.00$" to represent a positive number $< 0.005$, ``$0$" to represent zero, and we highlight in bold and with additional digits the relevant column that defines the state. If $\epsilon>8\cdot 10^{-5}$, then these two states are $\epsilon$-indistinguishable and could be merged together. In practice we observed that merging states such as the ones presented here has little impact in the quality of the solution, but substantially helps in controlling the size of the diagram. $\hfill \blacksquare$
\end{example}

In Example~\ref{ex:merging}, we saw a case with a tridiagonal matrix where, for $\epsilon$-exact diagrams, the value assignment made at layer $\ell=1$ could be irrelevant to computations in layer $\ell=7$. In \S\ref{sec:bandedMatrix} we formalize this intuition for problems with banded matrices. In particular, we show that for $\epsilon$-exact decision diagrams with $\epsilon>0$, the numbers of nodes at each layer can be bounded by a quantity that depends on the precision $\epsilon$, the condition number of $\bm{Q}$ and the bandwidth of the matrix, but is independent of the dimension $n$ of the problem. In other words, if the relevant parameters of the matrix $\bm{Q}$ are fixed, then the size of an $\epsilon$-exact decision diagrams is \emph{linear} in $n$. Leveraging these insights, we show that it is possible to obtain solutions to problem \eqref{eq:optRepeated} with objective at most OPT+$\varepsilon$, where OPT denotes the optimal objective, with runtime polynomial in $n$ and $1/\varepsilon$.

\section{A fully polynomial time approximation scheme}\label{sec:bandedMatrix}
In this section we derive FPTAS for \eqref{eq:miqo} where the bandwidth $k$ of $\bm{Q}$ and its condition number are fixed. Recall that the condition number of $\bm Q$ is denoted by $\cond(\bm Q)\triangleq \dfrac{\gamma_{\max}(\bm Q)}{\gamma_{\min}(\bm Q)}$, where $\gamma_{\max}(\bm Q)$ and $\gamma_{\min}(\bm Q)$ are the largest and smallest eigenvalues of $\bm Q$ respectively. 

\begin{remark}\label{rem:condition} Consider 
matrices of the form $\bm{Q}=\bm{I}+\bm{R}$ arising in monitoring problems \eqref{eq:inference}, where $\bm{R}$ encodes the temporal regularization. The minimum eigenvalue is $\gamma_{\min}(\bm{Q})=\min_{\bm{x}: \|\bm{x}\|_2^2=1}\|\bm{x}\|_2^2+\bm{x^\top R x}=1$
 if $R$ is given from a $k$-th order difference \eqref{eq:kth-Diff} with $k\geq 1$ or a moving average filter \eqref{eq:movingAverage}, where the constant $\bm{x^*}=\frac{1}{\sqrt{n}}\bm{1}$ is the associated optimal solution. The condition number is thus given by $$\cond(\bm{Q})=\gamma_{\max}(\bm{Q})=1+\lambda \gamma_{\max}(\bm{\bar R})$$ where $\lambda$ is the smoothness parameter and $\bm{\bar R}$ corresponds to the temporal filter (excluding $\lambda$) in \eqref{eq:kth-Diff} or \eqref{eq:movingAverage}. In practice, $\gamma_{\max}(\bm{\bar R})$ is typically small. For example, in the problems with $n=200$ with the moving average filter we use to generate Figure~\ref{fig:timesOnline}, we find $\gamma_{\max}(\bm{\bar R})=2.87$ if $k=2$ and $\gamma_{\max}(\bm{\bar R})=2.78$ if $k=3$; moreover, the dense moving average filter with $k=200$ yields $\gamma_{\max}(\bm{\bar R})=2.76$. Finally, note that a simple upper bound on the maximum eigenvalue when $\bm{\bar R}$ has bandwidth $k$ is $ \displaystyle\gamma_{\max} (\bm{\bar R})\le\norm{1}{\bm Q}\defeq\max_{i\in[n]}\sum_{j\in[n]}|Q_{ij}|\le k\max_{i,j\in [n]}|Q_{ij}|$, although this bound is rarely tight. $\hfill \blacksquare$
\end{remark}

Throughout this section we will compute upper bounds on the differences between states of compressed decision diagrams and approximations obtained by merging states (such as $\epsilon$-exact diagrams). Recall that the states store columns of special pseudoinverses associated with matrix $\bm{Q}$, corresponding to inverses of submatrices and padding remaining elements with zero. In order to use properties of inverses matrices, which are simpler than corresponding properties of pseudo-inverses, we introduce the following auxiliary definition.

\begin{definition}[padding matrix]
    For a $n$-dimensional matrix $\bm V\succeq 0$ and a subset $S\subseteq[n]$, the padding matrix ${\bm \hat V_S}\in \R^{n\times n}$ is defined as 
    \[ (\hat V_S)_{ij}=\begin{cases}
        V_{ij} & \text{if } i,j\in S\\
        1 & \text{if } i=j\notin S\\
        0&\text{otherwise}.
    \end{cases} \]
\end{definition}
\begin{example}
    Assume $n=5$, $k=2$, $S_1=\{1,2,4,5\}$, $S_2=\{1,4,5\}$ and $S_3=\{1,3,4\}$, then
    \[\footnotesize\bm V=\begin{pmatrix}
        *&*&*&&\\
        *&*&*&*&\\
        *&*&*&*&*\\
        &*&*&*&*\\
        &&*&*&*
    \end{pmatrix},\quad \bm{\hat V_{S_1}}=\begin{pmatrix}
        *&*&&&\\
        *&*&&*&\\
        &&1&&\\
        &*&&*&*\\
        &&&*&*
    \end{pmatrix},\quad \bm{\hat V_{S_2}}=\begin{pmatrix}
        *&&&&\\
        &1&&&\\
        &&1&&\\
        &&&*&*\\
        &&&*&*
    \end{pmatrix}, \quad \bm{\hat V_{S_3}}=\begin{pmatrix}
        *&&*&&\\
        &1&&&\\
        *&&*&*&\\
        &&*&*&\\
        &&&&1
    \end{pmatrix}. \hfill\blacksquare \]
\end{example}
 Observe that \emph{(i)} the bandwidth of the padding matrix does not increase; \emph{(ii)} given $\bm{\bar z}\in \{0,1\}^n$ and letting $S=\left\{i\in [n]:\bar z_i=1\right\}$, we find that $\left(\bm{ \hat V_S}\right)^{-1}_{ij}=\left(\bm{V}\circ \bm{\bar z\bar z^\top}\right)^{\dagger}_{ij}$ for all $i,j\in S$.

\subsection{Effects of one-shot merging}
Suppose we have constructed a compressed decision diagram up to layer $\ell$, and we proceed to merge $\epsilon$-indistinguishable nodes in this last layer. We now show that, after merging, the number of nodes in the last layer is bounded above by a quantity which does not depend on the layer $\ell$. The results rely on the following couple of lemmas. The first lemma, from the literature, states that the magnitude of the off-diagonal elements of the inverse of a banded matrix decay exponentially with the distance from the diagonal.

\begin{lemma}[\cite{demko1984decay}]\label{lem:decay banded inverse}
    Assume $\bm V\in\R^{n\times n}$ is a positive definite matrix with bandwidth $k$. Then the inequality  \[ \left| V^{-1}_{ij} \right|\le C_0\gamma^{\frac{|i-j|}{k}}\;\;\forall i,j\in [n] \]
    holds,
	where $\displaystyle \gamma=\frac{\sqrt{\cond(\bm V)}-1}{\sqrt{\cond(\bm V)}+1},\; C_0=\max\left\{ 1,\;\left( 1+\sqrt{\cond(\bm V)} \right)^2/(2\cond(\bm V)) \right\}/\gamma_{\min}(\bm V). $
\end{lemma}
Observe that since $\cond(\bm{V})\geq 1$, we find that $1/\gamma_{\min}(\bm{V})\leq C_0\leq 2/\gamma_{\min}(\bm{V})$. In monitoring problems, where the objective matrix satisfies $\gamma_{\min}(\bm{Q})=1$ (Remark~\ref{rem:condition}), we find that $C_0$ can be simply treated as a small constant independent of the condition number or other characteristics of the problem. 

Note that if $\bm V$ is a positive semidefinite matrix, then for all $i,j\in [n]$, $|V_{ij}|\le \VDiagMax\defeq \max_{h\in [n]} V_{hh}$. The second lemma establishes that, for a banded matrix $\bm{V}$, the difference between the last columns of $\bm{V^{-1}}$ and certain pseudoinverses associated with it is small.
\begin{lemma}\label{lem:bound diff of matrix inverses}
Let $\bm V\in\R^{n\times n}$ be a positive definite matrix with bandwidth $k$ partitioned as $\bm V=\begin{pmatrix}
    \bm A&\bm B^\top\\
    \bm B&\bm G
\end{pmatrix}$,
where $\bm{A}\in \R^{(n-m)\times (n-m)}$ and $\bm{G}\in \R^{m\times m}$ with $k<m<n$. Then, letting $S= \{n-m+1,\dots,n\}$, we have that 
\[ \left|V^{-1}_{ij}-\left(\bm{ \hat V_S}\right)^{-1}_{ij}\right|\le \begin{cases}
    C_0\gamma^{\frac{j-i}{k}}&\text{ if }1\le i\le n-m,\;j\ge n-m+1\\
    \VDiagMax C_0^2\left(\frac{1-\gamma^{\frac{n-m}{k}}}{1-\gamma^{\frac{1}{k}}}\right)^2\gamma^{\frac{i+j+2m-2n}{k}} &\text{ if }i\ge n-m+1,\;j\ge n-m+1.
\end{cases}\]

\end{lemma}
\proof
    By Lemma~\ref{lem:blockwiseInversion}, one has 
$\footnotesize
		\bm V^{-1}
=\begin{pmatrix}
	(\bm V/\bm G)^{-1}&\bm P\\ \bm P^\top&\bm G^{-1}+ \bm P^\top(\bm V/\bm G)\bm P
\end{pmatrix},
$
 where $\bm V/\bm G\defeq\bm A-\bm B^\top \bm{G^{-1}}\bm B$ is the Shur complement of the block $\bm G$ in $\bm V$ and $\bm{P}=-(\bm V/\bm G)^{-1}\bm B^{\top}\bm G^{-1}$. Moreover, define $\bm{H}\defeq \bm{ \hat V_S^{-1}}={\footnotesize\begin{pmatrix}
	\bm{I}&\bm 0\\ \bm 0&\bm{G^{-1}}
\end{pmatrix}.}$ By Lemma~\ref{lem:decay banded inverse}, one has $|P_{ij}|\le C_0\gamma^{\frac{j-i}{k}}$ for all $i\in [n-m]$ and $j\in [m]$, implying that the conclusion holds for $1\le i\le n-m$. If $i\ge n-m+1$ and $j\ge n-m+1$, one can deduce that
\begin{align*}
	\left|V^{-1}_{ij}-H_{ij}\right|=&\left|\sum_{j_1=1}^{n-m}\sum_{j_2=1}^{n-m}P_{j_1i}(\bm V/\bm G)_{j_1j_2}P_{j_2j}\right|\\
 \le & \VDiagMax\sum_{j_1=1}^{n-m}\sum_{j_2=1}^{n-m}\left|P_{j_1i}P_{j_2j}\right|\tag{$\because$ $|(\bm V/\bm G)_{j_1 j_2}|\le  A_{hh}\leq V_{\max}\;\forall h$}\\
	\le& \VDiagMax C_0^2\sum_{j_1=1}^{n-m}\gamma^{\frac{i-j_1}{k}}\sum_{j_2=1}^{n-m}\gamma^{\frac{j-j_2}{k}}\\
 =& \VDiagMax C_0^2\left(\gamma^{\frac{i+m-n}{k}}\frac{1-\gamma^{\frac{n-m}{k}}}{1-\gamma^{1/k}}\right)\cdot\left(\gamma^{\frac{j+m-n}{k}}\frac{1-\gamma^{\frac{n-m}{k}}}{1-\gamma^{1/k}}\right)\\
 =&\VDiagMax C_0^2\left(\frac{1-\gamma^{\frac{n-m}{k}}}{(1-\gamma^{1/k})}\right)^2\gamma^{\frac{i+j+2m-2n}{k}}.
\end{align*}
This finishes the proof. 
\endproof

We now show that the number of nodes obtained after merging $\epsilon$-indistinguishable nodes in an arbitrary layer of a decision diagram is bounded by a quantity independent of the layer. In particular, the proposition establishes that paths that coincide in the most recent $m$ value assignments, where $m$ is a sufficiently large number, can be merged into a single node in $\epsilon$-exact diagrams.
\begin{proposition}
Given a matrix $\bm{Q}$ with bandwidth $k$ and the associated compressed decision diagram $\dd$, there exists a constant $K$ (that depends only on the minimum and maximum eigenvalues of $\bm{Q}$ and bandwidth) such that for any
$m\geq k-1+k\frac{\log(K/\epsilon)}{|\log(\gamma)|}$,
any two paths representing solutions $\bm{\bar z^{\ell+1}}$ and $\bm{\hat z^{\ell+1}}$ such that $\hat z_i^\ell=0$ for all $i\leq\ell-m$ and $\bar z_i^\ell=\hat z_i^\ell$ for all $i>\ell-m$ lead to $\epsilon$-indistinguishable states.

\end{proposition}\label{prop:exactBound}
\proof{Proof} First assume that $\bm{\bar z^{\ell+1}}=\bm{1}$. Observe that since the matrix has bandwidth $k$, the relevant columns have indexes $j\geq \ell+1-k$. Then note that the non-zero rows of the state $\bm{\bar W^{\ell+1}}$ associated with $\bm{\bar z^{\ell+1}}$ are given by the relevant columns of the inverse of $\bm{V}=\bm{Q_{[\ell]}}$ (the submatrix of $\bm{Q}$ corresponding to the first $\ell$ indexes). Similarly, the non-zero rows of the state $\bm{\hat W^{\ell+1}}$ associated with $\bm{\hat z^{\ell+1}}$ are given by the relevant columns of the inverse of the padding matrix $\bm{\hat{V}_S}$ where $S=\{\ell-m+1,\dots,\ell\}$. Then, using Lemma~\ref{lem:bound diff of matrix inverses} we find that if $1\leq i\leq \ell-m$ and $j\geq \ell-m+1$, then for the relevant columns $C_0\gamma^{\frac{j-i}{k}}\leq C_0\gamma^{\frac{(\ell+1-k)-(\ell-m)}{k}}=C_0\gamma^{\frac{1+m-k}{k}}$. Similarly, if $i\geq \ell-m+1$ and $j\geq \ell-m+1$ then  \begin{align*}\VDiagMax C_0^2\left(\frac{1-\gamma^{\frac{\ell-m}{k}}}{1-\gamma^{\frac{1}{k}}}\right)^2\gamma^{\frac{i+j+2m-2n}{k}}\leq&  \VDiagMax C_0^2\frac{1}{\left(1-\gamma^{\frac{1}{k}}\right)^2}\gamma^{\frac{(\ell-m+1)+(\ell+1-k)+2m-2\ell}{k}}\\
=&\VDiagMax C_0^2\frac{\gamma^{\frac{1}{k}}}{\left(1-\gamma^{\frac{1}{k}}\right)^2}\gamma^{\frac{1+m-k}{k}}.\end{align*} Letting $K=\max\left\{C_0,\VDiagMax C_0^2\frac{\gamma^{\frac{1}{k}}}{\left(1-\gamma^{\frac{1}{k}}\right)^2}\right\}$, we find that the two nodes are $\epsilon$-indistinguishable when 
\begin{align*}
K\gamma^{\frac{1+m-k}{k}}\leq \epsilon\;\Leftrightarrow& \frac{1+m-k}{k}\log (\gamma)\leq \log(\epsilon/K)\\
\Leftrightarrow& \frac{1+m-k}{k}\geq \frac{\log(\epsilon/K)}{\log(\gamma)}\tag{$\because$ $\gamma<1\implies \log(\gamma)<0$}\\
\Leftrightarrow& m\geq k-1+k\frac{\log(\epsilon/K)}{\log(\gamma)}.
\end{align*}
Since $\gamma$ is a function of the eigenvalues of $\bm{Q}$ and $\VDiagMax\leq \gamma_{\max}(\bm{Q})$, the claim follows for this case.

To prove the other cases, let $S=\left\{i\in [\ell]:\hat z_i^{\ell+1}=1\right\}$, $T=\left\{i\in [\ell]:\bar z_i^{\ell+1}=1\right\}\supseteq S$ and let $\bm{U}=\bm{Q_{[\ell]}}$. The result can be obtained in an identical fashion using Lemma~\ref{lem:bound diff of matrix inverses} with $\bm{V}=\bm{\hat{U}_T}$ (corresponding to the state associated with $\bm{\bar z^{\ell+1}})$, since $\bm{\hat{V}_S}$ is the state associated with $\bm{\hat z^{\ell+1}}$.

\endproof

\begin{corollary}\label{cor:mergingSize}
Given a matrix $\bm{Q}$ with bandwidth $k$ and associated (exact) decision diagram $\dd$, a complete merging of $\epsilon$-indistinguishable states at layer $\ell+1>k$ results in at most $2^{k}(K/\epsilon)^{k/|\log(\gamma)|}$ nodes in this layer.
\end{corollary}
\proof
Merging first states that coincide in the most recent  $m\geq \lceil k-1+k\frac{\log(K/\epsilon)}{|\log(\gamma)|}\rceil $ assignments, where $K$ is the constant from Proposition~\ref{prop:exactBound}, we find that there are at most $2^{m}$ nodes in this layer, one for each possible combination of assignments in the last layers. 
\endproof

We point out, as noted by \citet{demko1984decay}, that the bounds in Lemma~\ref{lem:decay banded inverse} are not tight, thus the ensuing bounds are also not tight in general and $\epsilon$-exact merging will typically reduce the number of states by a much larger amount than the bound suggested in Corollary~\ref{cor:mergingSize}. Moreover, we naturally expect this bound to hold throughout all layers of an $\epsilon$-exact decision diagram and not just for a one-shot merging as proved above. We indeed observed this phenomenon in our computations, where the number of nodes per layer stabilizes in a quantity independent of $\ell$ and $n$. However, theoretically tracking the states merged in $\epsilon$-exact diagram is considerably challenging as it is highly dependent on the matrix $\bm{Q}$ of the problem. Nonetheless we show that, by altering the state of the decision diagram, we can design an algorithm with provable optimality guarantees.


\subsection{Truncated decision diagrams}

Motivated by the observation that $\epsilon$-exact decision diagrams merge states that coincide in the most recent $m$ value assignments, we define a new class of decision diagrams. We assume throughout this subsection that the problem is unconstrained, i.e., $Z=\{0,1\}^n$.

\begin{definition}[$m$-truncated decision diagram]\label{def:ddTruncated}
An $m$-truncated decision diagram is any decision diagram produced layer by layer according to Definition~\ref{def:fullDD} and \emph{merging} states that coincide in the most recent $m$ values assignments (selecting the state with zero value assignments in the first $\ell-1-m$ layers as the representative). 
\end{definition}

An $m$-truncated decision diagram has at most $2^m$ nodes per layer and thus at most $n2^{m+1}$ arcs overall. Theorem~\ref{theo:fptas} below states the main results of this section, namely that truncated diagrams can be used to design a FPTAS for problem \eqref{eq:optUnconstrained}, where matrix $\bm{Q}$ has fixed bandwidth and minimum/maximum eigenvalues.

\begin{theorem}\label{theo:fptas}
Let $h^*$ be the optimal objective value of problem \eqref{eq:optRepeated}. There exists a constant $C$ that depends only on the bandwidth $k$ and the minimum/maximum eigenvalue of $\bm{Q}$ such that the shortest path between the root and a terminal in
    an $m$-truncated decision diagram with $\displaystyle m=\left\lceil\frac{k}{|\log(\gamma)|}\log\left( \frac{C\|\bm{d}\|_\infty^2 n}{\varepsilon} \right)\right\rceil$ has length $h_m^*$ such that $|h_m^*-h^*|\leq \varepsilon$. 
\end{theorem} 
Since the number of arcs in the diagram is bounded by $\displaystyle n\left\lceil \left( \frac{C\|\bm{d}\|_\infty^2 n}{\varepsilon} \right)^{1+\frac{k}{|\log(\gamma)|}}\right\rceil$, we find that a shortest path in this problem is indeed polynomial in the dimension $n$ and the precision $\|\bm{d}\|_\infty^2/\varepsilon$ whenever the bandwidth $k$ and minimum/maximum eigenvalues of $\bm{Q}$ are fixed.

From the definition of truncated decision diagrams, we find that Lemma~\ref{lem:bound diff of matrix inverses} can be used to bound the differences between states in the complete decision diagram $\dd_{\text{full}}$ and $m$-truncated diagram $\dd_{m}$ reached by paths representing the same solutions. We formally state this fact in Proposition~\ref{prop:diff-W}.

\begin{proposition}\label{prop:diff-W}
Given any partial solution $\bm{\bar z^\ell}\in \{0,1\}^n$ represented by path $(a_1,a_2,\dots,a_{\ell-1})$ in a complete diagram $\dd_{\text{full}}$ and path $(\alpha_1,\alpha_2,\dots,\alpha_{\ell-1})$ in an $m$-truncated diagram $\dd_m$, and letting $\bm{\bar W^\ell}$ and $\bm{\bar \Omega^\ell}$ be the states reached by those paths in the complete and truncated diagrams respectively, the bound 
\begin{equation}\label{eq:comparison of W}
    \left| \bar{\Omega}^{\ell}_{ij} - \bar{W}^\ell_{ij} \right|\le\begin{cases}
    C_0\gamma^{\frac{j-i}{k}}&\text{ if }1\le i\le \ell-m-1,\;j\ge \ell-m\\
    C_1\gamma^{\frac{i+j+2m-2\ell+2}{k}} &\text{ if }i\ge \ell-m,\;j\ge \ell-m
\end{cases}
\end{equation}
holds, where $C_1=\QDiagMax \left(\frac{C_0}{1-\gamma^{\frac{1}{k}}}\right)^2$.
\end{proposition}

\ignore{
\proof
Denote $\alpha=\{i\in[n]: \bar z^\ell_i=1\}$ and $\alpha_m=\{i\in[n]: \Bar{z}_i^{\ell,m}=1 \}.$ If $i\notin\alpha$ or $j\notin\alpha$, then $\Bar{W}^{\ell, m}_{ij} = \Bar{W}^\ell_{ij}=0$. Thus, we assume $i\in\alpha$ and $j\in\alpha$.  By the definition of 1-padding matrices, it suffices to prove
\begin{equation}\label{eq:comparision of W aux}
    \left|  \left(\hat{\bm Q}^{\alpha_m}\right)^{-1}_{ij}-\left(\hat{\bm Q}^{\alpha}\right)^{-1}_{ij} \right|\le\begin{cases}
    C_0\gamma^{\frac{j-i}{k}}&\text{ if }1\le i\le \ell-m-1,\;j\ge \ell-m\\
    C_1\gamma^{\frac{i+j+2m-2\ell+2}{k}} &\text{ if }i\ge \ell-m,\;j\ge \ell-m.
\end{cases}
\end{equation}
 Letting $S=[\ell-1]\supseteq \alpha$ and $T=[\ell-1]\backslash[\ell-1-m]\supseteq \alpha_m$, one can deduce that
 $\left(\hat{\bm Q}^{\alpha_m}\right)_S = \begin{pmatrix}
     \bm I&\bm 0\\
     \bm 0& \left(\hat{\bm Q}^{\alpha}\right)_{T}
 \end{pmatrix}$. The conclusion \eqref{eq:comparision of W aux} follows from Lemma~\ref{lem:bound diff of matrix inverses}.

\endproof
}
Observe that in the proposition, we used the bound $1-\gamma^{\frac{\ell-m}{k}}\leq 1$, which is a good approximation if $\ell$ is large.
We now prove that not only the states stored in the complete and truncated diagrams are similar, but also the arc lengths computed once the cost vectors $\bm c, \bm d\in\R^n$ are realized. Recall that the arc lengths are computed as $l_a=c_{t_a}\nu_{t_a} - \frac{1}{2}(\bm{d^\top u_{a}})^2$ where $a$ is an arc in the complete diagram. Similarly, given an arc $\alpha$ in the truncated diagram, we let $l_\alpha$ denote its length in this diagram, which is computed identically. Proposition~\ref{prop:comparision of arc length} states that the differences in the lengths of certain arcs can be bounded.
\begin{proposition}\label{prop:comparision of arc length}
    Given any solution $\bar{\bm z}\in\{0,1\}^n$, let $(a_1,\dots, a_{n})$  and $(\alpha_1,\dots, \alpha_n)$ be the paths representing this solution in the complete diagram $\dd_{\text{full}}$ and the $m$-truncated diagram $\dd_m$, respectively, and define lengths according to Proposition~\ref{prop:pathLength}. Then, given any layer $\ell\in [n]$ 
    the inequality
    \begin{equation*}
        \left|l_{a_\ell}-l_{\alpha_\ell}\right| \le C \|\bm{d}\|_\infty^2 \gamma^{m/k},
    \end{equation*}
    holds, where $\displaystyle C\defeq \frac{C_0Q_{\max}\left( 2C_0+2C_1+ C_0C_1Q_{\max}\right)\gamma^{1/k}}{2(1-\gamma^{1/k})^2}$. 
\end{proposition}

\proof
 Let $\bm{ \bar W^\ell}$ and $\bm{ \bar \Omega^{\ell}}$ be the matrix stored by the state $h_{a_{\ell-1}}$ and $h_{\alpha_{\ell-1}}$ in $\mathcal{D}$ and $\mathcal{D}_m$, respectively. The conclusion trivially holds if either $\bar z_\ell=0$ since $|l_{a_\ell}-l_{\alpha_\ell}|=0$ or  $\ell\le m+1$ since no truncation is performed in this case. Thus, we assume $\ell\ge m+2$ and $\bar z_\ell=1$. 
By definition, one has
\begin{align}
    \left| l_{a_\ell} -l_{\alpha_\ell} \right|=&\frac{1}{2}\left| \left( \bm {d^\top u_{\alpha_\ell}} \right)^2 - \left( \bm{d^\top u_{a}} \right)^2\right|  = \frac{1}{2} \left| \bm{d^\top}\left( \bm {u_{\alpha_\ell}} +\bm{ u_{a_\ell}}\right)  \cdot  \bm{d^\top} \left(\bm{u_{\alpha_\ell}}-\bm{u_{a_\ell}}\right)\right|. \label{eq:expand-diff-l}
    \end{align}
In the following, we bound each multiplier in \eqref{eq:expand-diff-l}. First observe that from Remark~\ref{rem:uDiff} we find that $\bm{u_{a_\ell}}=\bm{\bar W_{\ell}^{\ell+1}}/\sqrt{\Bar{W}^{\ell+1}_{\ell\ell}}$, where $\bm{\bar W^{\ell+1}}\defeq\left(\bm Q\circ \bm{\hat{z}}\bm{\hat{z}}^\top\right)^{\dagger}$ and $\bm{\hat z}\defeq\sum_{i=1}^\ell\bar z_i\bm{e_i}$. Moreover, letting $\bm{\hat z^{m+1}}\defeq \sum_{i=\ell-m}^\ell \bar z_i\bm{e_i}$ and $ \bm{\bar\Omega^{\ell+1}}\defeq \left(\bm Q\circ \bm{\hat{z}^{m+1}}\left(\bm{\hat{z}^{m+1}}\right)^\top\right)^{\dagger}$, we find that $\bm{\bar\Omega^{\ell+1}} = \bm{\bar\Omega^{\ell}}+\bm{u_{\alpha_\ell} u_{\alpha_\ell}^\top}$ and $\bm{u_{\alpha_{\ell}}}=\bm{ \bar \Omega_{\ell}^{\ell+1}}/\sqrt{\bar{\Omega}_{\ell\ell}^{\ell+1}}$. In addition, we note that  $ 1/\QDiagMax\le \bar{W}_{\ell\ell}^{\ell+1}, \bar{\Omega}_{\ell\ell}^{\ell+1} \le 1/\gamma_{\min}(\bm Q)$.
\begin{itemize}
    \item Define $S=\{i\in[\ell]: \bar z_i=1\}$ and $S_{m+1}=\left\{i:\hat{z}^{m+1}_i=1\right\}$. Note that $S_{m+1}\subseteq S$. Since \begin{align*}
        \left|\left(\bm{u_{a_{\ell}}} \right)_i\right|=\left| \Bar{W}^{\ell+1}_{i\ell} \right|/\sqrt{\bar{W}_{\ell\ell}^{\ell+1}}=\left|\left(\bm{ \hat Q_S}\right)^{-1}_{i\ell}\right|/\sqrt{\bar{W}_{\ell\ell}^{\ell+1}}\le \frac{C_0}{\sqrt{\bar{W}_{\ell\ell}^{\ell+1}}}\gamma^{\frac{\ell -i}{k}}\quad\forall i\in[\ell], \tag{$\because$ Lemma~\ref{lem:decay banded inverse}} 
    \end{align*} we deduce that
    \begin{align*}
        \left|  \bm {d^\top \bm u_{a_\ell}} \right|=\left|\sum_{i=1}^{\ell} d_i\left(\bm {u_{a_{\ell}}} \right)_i \right|
        \le \|\bm{d}\|_\infty\frac{C_0}{\sqrt{\bar{W}_{\ell\ell}^{\ell+1}}}\sum_{i=1}^{\ell} \gamma^{\frac{\ell-i}{k}}\le \frac{C_0\|\bm{d}\|_\infty\sqrt{\QDiagMax}}{1-\gamma^{1/k}}.\label{eq:bound-begin}
    \end{align*}   
    With identical arguments applied to $S_{m+1}$, we can prove that $\displaystyle\sum_{i=1}^\ell\left|\bar \Omega_{i\ell}^{\ell+1}\right|\le \frac{C_0}{1-\gamma^{1/k}}$ and $\displaystyle
        \left| \bm{ d^\top u_{\alpha_\ell}} \right|  \le \frac{C_0\|\bm{d}\|_\infty \sqrt{\QDiagMax}}{1-\gamma^{1/k}}$.
    Thus, 
    \begin{equation}
        \left| \bm{d^\top}\left( \bm{u_{\alpha_\ell}} +\bm{u_{a_\ell}}\right)  \right| \le \left|  \bm{d^\top u_{a_\ell}} \right| + \left| \bm{d^\top u_{\alpha_\ell}} \right| \le \frac{2C_0\|\bm{d}\|_\infty\sqrt{\QDiagMax}}{1-\gamma^{1/k}}. \label{eq:bound-innerProd-sum}
    \end{equation}
    \item To bound $\left| \bm {d^\top} \left(\bm{\bm u_{a_\ell}}-\bm{u_{\alpha_\ell}}\right)\right|$, we find that
    \begin{align*}
        &\left| \bm {d^\top} \left(\bm{u_{\alpha_\ell}}-\bm{u_{a_\ell}}\right)\right| = \left|\sum_{i=1}^{\ell}d_i\left( \bar{W}^{\ell+1}_{i\ell}/\sqrt{\bar{W}_{\ell\ell}^{\ell+1}} - \bar{\Omega}^{\ell+1}_{i\ell}/\sqrt{\bar{\Omega}_{\ell\ell}^{\ell+1}} \right)\right|\\
        \le & \|\bm{d}\|_\infty\sum_{i=1}^\ell\left| \left( \bar{W}^{\ell+1}_{i\ell} -\bar{\Omega}^{\ell+1}_{i\ell} \right)/\sqrt{\bar{W}_{\ell\ell}^{\ell+1}} + \bar{\Omega}^{\ell+1}_{i\ell}\left(1/\sqrt{\bar{W}_{\ell\ell}^{\ell+1}} - 1/\sqrt{\bar{\Omega}_{\ell\ell}^{\ell+1}}\right) \right|\\
        \le &\|\bm{d}\|_\infty\left[ \sum_{i=1}^{\ell}\left| \bar{W}^{\ell+1}_{i\ell} - \bar{\Omega}^{\ell+1}_{i\ell} \right|/\sqrt{\bar W_{\ell\ell}^{\ell+1}} +  \frac{\left|\sqrt{\bar{W}_{\ell\ell}^{\ell+1}} - \sqrt{\bar{\Omega}_{\ell\ell}^{\ell+1}}\right|}{\sqrt{\bar{W}_{\ell\ell}^{\ell+1}}\sqrt{\bar{\Omega}_{\ell\ell}^{\ell+1}} }  \sum_{i=1}^\ell\left|\bar \Omega_{i\ell}^{\ell+1}\right|\right].\end{align*}
Observe that (from the first part of the proof) we find that $\sum_{i=1}^\ell\left|\bar \Omega_{i\ell}^{\ell+1}\right|\le \frac{C_0}{1-\gamma^{1/k}}$. Moreover, $1/\sqrt{\bar W_{\ell\ell}^{\ell+1}}\leq \sqrt{Q_{\max}}$ and $1/\left(\sqrt{\bar{W}_{\ell\ell}^{\ell+1}}\sqrt{\bar{\Omega}_{\ell\ell}^{\ell+1}}\right)\leq Q_{\max}$. In addition, $\sqrt{\bar{W}_{\ell\ell}^{\ell+1}} - \sqrt{\bar{\Omega}_{\ell\ell}^{\ell+1}}=\frac{\bar{W}_{\ell\ell}^{\ell+1} - \bar{\Omega}_{\ell\ell}^{\ell+1}}{\sqrt{\bar{W}_{\ell\ell}^{\ell+1}} + \sqrt{\bar{\Omega}_{\ell\ell}^{\ell+1}}}$. Using these substitutions and bounds, we find that $\left| \bm {d^\top} \left(\bm{u_{\alpha_\ell}}-\bm{u_{a_\ell}}\right)\right|$ is less or equal than
        \begin{align*}
        &\|\bm{d}\|_\infty\left[ \sqrt{\QDiagMax}\left(\sum_{i=1}^{\ell -m-1} \left| \bar W^{\ell+1}_{i\ell}-\bar \Omega_{i\ell}^{\ell+1} \right|+\sum_{i=\ell-m}^{\ell}\left| \bar W^{\ell+1}_{i\ell}-\bar \Omega_{i\ell}^{\ell+1} \right|\right) + \frac{C_0Q_{\max}}{1-\gamma^{1/k}} \frac{\left|\bar{W}_{\ell\ell}^{\ell+1} - \bar{\Omega}_{\ell\ell}^{\ell+1}\right|}{\sqrt{\bar{W}_{\ell\ell}^{\ell+1}} + \sqrt{\bar{\Omega}_{\ell\ell}^{\ell+1}}}  \right]\\
        \le & \|\bm{d}\|_\infty\left[\sqrt{Q_{\max}} \left(\sum_{i=1}^{\ell -m-1}\! C_0\gamma^{\frac{\ell-i}{k}} + \!\!\!\! \sum_{i=\ell-m}^{\ell}\! C_1\gamma^{\frac{i+2m-\ell+2}{k}}\right)+\frac{C_0Q_{\max}^{3/2}}{2(1-\gamma^{1/k})}C_1\gamma^{\frac{2m+2}{k}}\right] \tag{Proposition~\ref{prop:diff-W}}\\
        \le &    \|\bm{d}\|_\infty\left[\sqrt{Q_{\max}} \left( C_0\frac{\gamma^{\frac{m+1}{k}}}{1-\gamma^{1/k}} + C_1\frac{\gamma^{\frac{m+2}{k}}}{1-\gamma^{1/k}}\right)+\frac{C_0Q_{\max}^{3/2}}{2(1-\gamma^{1/k})} C_1\gamma^{\frac{2m+2}{k}}\right]\\
        \le & \frac{\|\bm{d}\|_\infty\sqrt{Q_{\max}}\left( 2C_0+2C_1+ C_0C_1Q_{\max}\right)}{2(1-\gamma^{1/k})} \gamma^{\frac{m+1}{k}}, \numberthis \label{eq:bound-innerProd-diff}
    \end{align*}
\end{itemize}
where the last inequality is due to $\gamma\le1$. Plugging the bounds \eqref{eq:bound-innerProd-sum} and \eqref{eq:bound-innerProd-diff} in \eqref{eq:expand-diff-l}, we conclude that
\begin{align*}
    \left| l_{a_\ell} -l_{\alpha_\ell} \right| \le \frac{1}{2} \frac{2C_0\|\bm{d}\|_\infty\sqrt{\QDiagMax}}{1-\gamma^{1/k}} \frac{\|\bm{d}\|_\infty\sqrt{Q_{\max}}\left( 2C_0+2C_1+ C_0C_1Q_{\max}\right)}{2(1-\gamma^{1/k})} \gamma^{\frac{m+1}{k}} = C\|\bm{d}\|_\infty^2\gamma^{m/k},
\end{align*}
finishing the proof.

\endproof

Since the lengths of corresponding arcs in both diagrams are similar, we can conclude that the shortest in both diagrams also has similar length.
\begin{corollary}\label{prop:comparison of path length}
    Let $h^*$ and $h_m^*$ be the lengths of the shortest paths between roots and a terminal node in the full diagram $\dd_{\text{full}}$ and truncated diagram $\dd_m$, respectively, with arc lengths defined according to Proposition~\ref{prop:pathLength}. Then the inequality
    $ |h^*-h^*_m|\le C\|\bm{d}\|_\infty^2 n\gamma^{m/k} $ holds.
\end{corollary}
\proof
Assume $(a_1,\dots, a_n)$ is a shortest path in $\dd_{\text{full}}$, and let $(\alpha_1,\dots, \alpha_n)$ be the path in $\mathcal{D}_m$ representing the same solution. Then
    $h^* = \sum_{\ell=1}^n l_{a_{\ell}}$ and $h^*_m\le \sum_{\ell=1}^n l_{\alpha_{\ell}}$,
implying
\[ h^*_m-h_*\le \sum_{\ell=1}^n \left(l_{\alpha_{\ell}} - l_{a_{\ell}}\right) \le \sum_{\ell=1}^n C\|\bm{d}\|_\infty^2 \gamma^{m/k}= C_2\|\bm{d}\|_\infty^2n\gamma^{m/k},\]
where the last inequality is due to Proposition~\ref{prop:comparision of arc length}. Using the same argument, one can prove the other direction $h_*-h_*^m\le C\|\bm{d}\|_\infty^2n\gamma^{m/k}$.

\endproof
 We can now conclude the proof of Theorem~\ref{theo:fptas}.

\proof{( Theorem~\ref{theo:fptas})}
Following from Corollary~\ref{prop:comparison of path length}, it is sufficient to set $m$ such that $C\|\bm{d}\|_\infty^2n\gamma^{m/k}\le\varepsilon$, leading to the result.

\endproof

\section{Convexification}\label{sec:convexification}

So far, the main focus has been in solving problem \eqref{eq:miqo}. In this section, we show that the construction of the decision diagram is equivalent to the convexification of the constraints of the problem, that is, the convexification of set
$$X_{Q,Z}\defeq\left\{(\bm{x},\bm{z},x_0)\in \R^n\times \{0,1\}^n\times \R: x_0\geq \bm{x^\top Qx},\; \bm{x}\circ (\bm{1}-\bm{z})=\bm{0},\; \bm{z}\in Z\right\}$$
in an extended space. Let $\dd=(\mathcal{G},\nu,\bm{u})$ with $\mathcal{G}=(N,A)$ be the compressed decision diagram obtained according to Definition~\ref{def:compressedDD} with a potentially enlarged state space to include constraints $z\in Z$, as discussed in Remark~\ref{rem:constraints}.  Recall that given an arc $a\in A$, $h_a,t_a\in N$ denote the head and tail of the arc and $\ell(a)$ denotes the layer from which the arc emanates, thus $\bm{e_{\ell(a)}}\in \{0,1\}^n$ is the vector that has a one in position $\ell(a)$ and $0$ elsewhere. Given a node $v\in N$, we denote by $\ell(v)$ the layer in which the node is. Moreover, in this section we adopt the convention that $0/0=0$ and $c/0=+\infty$ if $c>0$. The next theorem states that the closure of the convex hull of $X_{Q,Z}$ admits an extended SOCP-representable formulation whose number of additional variables is linear in the number of arcs in the decision diagram.

\begin{theorem}\label{theo:hullDD}Point $(\bm{x},\bm{z},x_0)\in \conv(X_{Q,Z})$ if and only if there exist $\bm{r},\bm{w}\in \R^A$ such that the system
\begin{subequations}\label{eq:convexhull}
\begin{align}
    &x_0\geq \sum_{a\in A}\frac{w_a^2}{r_a}\label{eq:convexhull_obj}\\
    &\bm{x}=\sum_{a\in A}\bm{u_a}w_a\label{eq:convexhull_inverse}\\
    &\bm{z}=\sum_{a\in A:\nu_a=1}\bm{e_{\ell(a)}}r_a\label{eq:convexhull_assignment}\\
    &\sum_{a\in A: h_a=v}r_a=\sum_{a\in A: t_a=v}r_a\qquad \forall v\in N:2 \leq \ell(v)\leq n\label{eq:convexhull_f3}\\
    &\sum_{a\in A:\ell(a)=1}r_a = 1; \; \sum_{a\in A:\ell(a)=n}r_a = 1; \; \bm{r}\geq \bm{0}\label{eq:convexhull_nonneg}
\end{align}
\end{subequations}
has a solution.
\end{theorem}

Intuitively, variables $\bm{r}$ in \eqref{eq:convexhull} indicate whether an arc in the decision diagram is used by a given solution.
Conditions \eqref{eq:convexhull_f3}-\eqref{eq:convexhull_nonneg} state that variable $\bm{r}$ is in the convex hull of all paths between the root in layer $\ell=1$ and a node in the last layer $\ell=n+1$; conditions \eqref{eq:convexhull_assignment} link paths and the indicator variables represented by the paths; terms $w_a^2/r_a$ in \eqref{eq:convexhull_obj} explicitly enforce the logical condition ``$r_a=0\implies w_a=0$", and implicitly enforce the indicator constraints ``$z_i=0\implies x_i=0$" as well; finally, \eqref{eq:convexhull_obj}-\eqref{eq:convexhull_inverse} ensure that $x_0\geq \bm{x^\top Qx}$ is satisfied. 

The rest of this section is devoted to the proof of Theorem~\ref{theo:hullDD}. The proof relies on the following recent result.
\begin{lemma}[\citet{wei2023convex}] \label{theo:convexification}If $\bm{Q}\succ 0$, then 
\begin{align*}
\conv(X_{Q,Z})=\Big\{(\bm{x},\bm{z},x_0)\in \R^{2n+1}: \exists \bm{W}\in \R^{n\times n} \text{such that }\begin{pmatrix}x_0 &\bm{x^\top}\\
\bm{x}&\bm{W}\end{pmatrix}\succeq 0,\; (\bm{z},\bm{W})\in \conv(P_{Q,Z})\Big\} 
\end{align*}
where 
$P_{Q,Z}\defeq\left\{\bm{z}\in Z,\,\bm{W}\in \R^{n\times n}:\bm{W}=\left(\bm{Q}\circ \bm{zz^\top}\right)^\dagger\right\}.$ 
\end{lemma}
Lemma~\ref{theo:convexification} indicates that the characterization of the non-polyhedral set $\cl\conv(X_{Q,Z})$ reduces to describing the polytope $\conv(P_{Q,Z})$ (defined as the convex hull of a finite set of points). The first step of our proof is thus to show that we can obtain an explicit description of $\conv(P_{Q,Z})$ from $\dd$.

\begin{proposition}\label{prop:hullPolytope}
A point $(\bm{z},\bm{W})\in \conv(P_{Q,Z})$ if and only if there exists $\bm{r}\in \R^A$ satisfying \eqref{eq:convexhull_assignment}-\eqref{eq:convexhull_nonneg} such that 
$\bm{W}=\sum_{a\in A}\bm{u_au_a^\top}r_a.$
\end{proposition}
\proof{(Proposition~\ref{prop:hullPolytope})}
Note that $$P_{Q,Z}=\left\{(\bm{z},\bm{W}):\exists \bm{r}\in \{0,1\}^A \text{ s.t. }\bm{W}=\sum_{i\in A}\left(\bm{u_au_a^\top}\right)r_a,\; \eqref{eq:convexhull_assignment}-\eqref{eq:convexhull_nonneg}\right\}.$$ 
Indeed, there is a one to one correspondence between points $\bm{z}\in Z$ and paths in $\mathcal{G}$. Moreover, any path in $\mathcal{G}$ corresponds to a binary point satisfying \eqref{eq:convexhull_f3}-\eqref{eq:convexhull_nonneg}. Constraints \eqref{eq:convexhull_assignment} force $\bm{z}$ to be the solution represented by the path, and the condition $\sum_{a\in A}\sum_{a\in A}\bm{u_au_a^\top}r_a=\left(\bm{W}\circ \bm{zz^\top}\right)^\dagger$ follows from the second statement in Proposition~\ref{prop:correctReduced} and cancelling out the terms in the telescoping sum induced by vectors $\bm{u}$ in the path. 

Now observe that the extreme points of the polytope $\bar P =\left\{(\bm{z},\bm{W},\bm{r}):\bm{W}=\sum_{i\in A}\left(\bm{u_au_a^\top}\right)r_a,\; \eqref{eq:convexhull_assignment}-\eqref{eq:convexhull_nonneg}\right\}$ are precisely the binary points (in $\bm{r}$) satisfying the constraints. Indeed, any extreme point is the unique optimal solution of an optimization problem of the form 
\begin{align*}
\min_{(\bm{z},\bm{W},\bm{r})\in \bar P}\bm{\alpha^\top x}+\bm{\beta^\top r}+\sum_{i=1}^n\sum_{j=1}^n \Gamma_{ij}W_{ij},
\end{align*} 
where $(\bm{\alpha},\bm{\beta}, \bm{\Gamma})$ are properly chosen. Using the equality constraints to project out $(\bm{z},\bm{W})$ we can rewrite the problem as linear optimization over variables $\bm{r}$ over constraints \eqref{eq:convexhull_f3}-\eqref{eq:convexhull_nonneg}: since those constraints are totally unimodular, the optimal solution is binary in $\bm{r}$. Since $\text{proj}_{\bm{r}}\bar P$ is a relaxation of $P_{Q,Z}$, and the extreme points of $\bar P$ are precisely the points in $P_{Q,Z}$, it follows that $\conv(P_{Q,Z})=\text{proj}_{\bm{r}}\bar P$.
\endproof

Equipped with Proposition~\ref{prop:hullPolytope}, we can conclude the proof of the theorem.
\begin{proof}{Proof of Theorem~\ref{theo:hullDD}}
From Lemma~\ref{theo:convexification} and Proposition~\ref{prop:hullPolytope}, we find that an extended formulation for $\text{cl conv}(X_{Q,Z})$ is given by 
\begin{align}\label{eq:hullSdp}
&\begin{pmatrix}x_0 &\bm{x^\top}\\
\bm{x}&\bm{W}\end{pmatrix}\succeq 0,\; 
\bm{W}=\sum_{a\in A}\bm{u_au_a^\top}r_a,\; \eqref{eq:convexhull_assignment}-\eqref{eq:convexhull_nonneg}.
\end{align}
Since $\bm{W}$ is defined via an equality constraint as the sum of given rank-one matrices multiplied by nonnegative variables, we immediately find \cite[][pp 227-229]{nesterov1994interior} that an SOCP representation of the set is \eqref{eq:convexhull}. 
\end{proof}

\section{Computations}\label{sec:computations}

In this section we present computational results with financial time series data. In \S\ref{sec:comp_instances} we present the instances used, in \S\ref{sec:comp_methods} we discuss the methods used, and in \S\ref{sec:comp_resOffline} and \S\ref{sec:comp_resOnline} we present the results in offline and online settings, respectively. We did not observe any major numerical difficulties, and report in Appendix~\ref{sec:numPrec} results concerning the numerical precision of methods tested.

\subsection{Instances}\label{sec:comp_instances}

Given a signal $\bm{y}\in \R^n$, smoothness parameter $\lambda\in \R_+$, sparsity parameter $\mu\in \R_+$, width parameter $k\in \Z_+$ and contiguity parameter $\tau\in \Z_+$, we consider inference problems with a moving average filter and constraints on the number of contiguous nonzero values as described by constraints \eqref{eq:consecutiveOnes}, that is, problems of the form
\begin{subequations}\label{eq:denoising}
\begin{align}
\min_{\bm{x},\bm{z},\bm{\zeta}}\;& \sum_{i=1}^n (y_i-x_i)^2 +\lambda \sum_{i=2}^n \left(x_i-\frac{1}{\min\{k,i-1\}}\sum_{j=1}^{\min\{k,i-1\}}x_{i-j}\right)^2+\mu\sum_{i=1}^n z_i\\
\text{s.t.}\; & \eqref{eq:consecutiveFirst}-\eqref{eq:consecutiveLast} \\
&\bm{x}\in \R^n,\;\bm{z}\in \{0,1\}^n,\;\bm{\zeta}\in \{0,1\}^{n+1-\tau}.
\end{align}
\end{subequations}
By convention we let \eqref{eq:denoising} with $\tau=0$ be the problem with no contiguity constraints, that is, obtained by removing constraints \eqref{eq:consecutiveFirst}-\eqref{eq:consecutiveLast} and variables $\bm{\zeta}$. The signal $\bm{y}$ are obtained from financial data, as described next. Time periods where $x_i\neq 0$ and $z_i=1$ in optimal solutions of \eqref{eq:denoising} correspond to high-volatility periods, which correlate with economic downturns \citep{campajola2022modelling}. 

\subsubsection*{Dataset} We use the financial data provided by Boris Marjanovic in Kaggle \footnote{\url{https://www.kaggle.com/datasets/borismarjanovic/price-volume-data-for-all-us-stocks-etfs}} to construct inputs $\bm{y}$. The original data contains the daily price for all US-based stocks and ETFs trading on the NYSE, NASDAQ, and NYSE MKT from 1990 and 2017. We then process the data as follows:
\begin{itemize}
    \item We only keep the stocks and ETFs who traded in at least 99\% of the days reported. For missing days, we set the price to be the same as the price for the most recently reported day.
    \item We compute the daily changes in the stocks and ETFs. The resulting dataset, containing 214 securities and 7,022 time periods, can be found online at \url{https://sites.google.com/usc.edu/gomez/data}.
\end{itemize}

From this data, we construct the instances as follows, depending on the setting.

\paragraph{Offline setting}
In this setting, the goal is to solve problem \eqref{eq:denoising} once. Given the dimension parameter $n$ and a security of interest, we create a signal $\bm{y}\in \R^n$ by partitioning the time horizon in $n$ epochs corresponding each to $\lfloor 7,022/n\rfloor$ days (the last $7,022-\lfloor 7,022/n\rfloor n$ observations are discarded), and $y_i$ is set to be the average daily change of the chosen security in the $i$-th epoch. Finally, data are standardized so that $\|\bm{y}\|_2^2=1$ and $\sum_{i=1}^n y_i=0$.

\paragraph{Online setting}
In this setting, we use the disaggregated signal $\bm{y_\textbf{full}}\in \R^{7022}$, standardized so that $\|\bm{y_\textbf{full}}\|_2^2=1$ and $\sum_{i=1}^{7022} (y_\text{full})_i=0$. Then, letting $\bm{y^t}\in \R^{200}$ denote the $200$-dimensional subvector of $\bm{y_\textbf{full}}$ starting at position $t$, that is, $y_{1+i}^t=(y_\text{full})_{t+i}$ for $0\leq i\leq 199$, the goal is to sequentially solve \eqref{eq:denoising} to optimality for all signals $\bm{y^t}$, $t=1,\dots, 6823$. This experiment emulates the situation where data is made available one datapoint at the time, and a decision-maker resolves the filtering problem \eqref{eq:denoising} on the fly to quickly detect high-volatility moments. 

\subsection{Methods and metrics}\label{sec:comp_methods}

We compare the following methods for offline settings.

\noindent $\bullet$ \textbf{Mosek} Mosek 10.0 branch-and-bound solver, formulating \eqref{eq:denoising} using the perspective reformulation as 
\begin{subequations}
\begin{align*}
\|\bm{y}\|_2^2+\min_{\substack{\bm{x}\in \R^n,\bm{z}\in \{0,1\}^n\\\bm{\zeta}\in \{0,1\}^n \bm{t}\in \R_+^n}}\;& \sum_{i=1}^n (-2y_ix_i+t_i) +\lambda \sum_{i=2}^n \left(x_i-\frac{1}{\min\{k,i-1\}}\sum_{j=1}^{\min\{k,i-1\}}x_{i-j}\right)^2+\mu\sum_{i=1}^n z_i\\
\text{s.t.}\;&x_i^2\leq t_iz_i\qquad \forall i\in [n]\\
&\eqref{eq:consecutiveFirst}-\eqref{eq:consecutiveLast}.
\end{align*}
\end{subequations}
The number of threads is set to one and a time limit of 1,800 seconds is imposed, and the default settings are used otherwise. For this method we compute the solution \underline{time} in seconds (timeouts count as 1,800 seconds), the number of branch-and-bound \underline{nodes} explored until proving optimality or the time limit is reached, and the \underline{percentage} of instances solved to optimality. 

\noindent $\bullet$ \textbf{Decision diagram} Solving problem \eqref{eq:denoising} by constructing a $10^{-5}$-exact decision diagrams according to Definition~\ref{def:ddEpsilon}. We also set a time limit of 1,800 seconds, and abort the construction of the diagram if the time limit is reached. We report the number of arcs \underline{$|A|$} of the diagram, the time in seconds used to construct the decision diagram (\underline{time\_dd}) --time limits count as 1,800 seconds--, the time in seconds required to solve the problem as a shortest 
path problem after construction of the diagram(\underline{time\_sp}), and the \underline{percentage} of instances solved; finally, we also report the time required to solve the SOCP relaxations ensuing from Theorem~\ref{theo:hullDD} using Mosek (\underline{time\_hull}).

For online settings, where each run requires the solution of 6,823 MIO problems to optimality, it is not realistic to use Mosek. Thus, we only report the times associated with constructing the decision diagram once (which can be done offline) and then solving each online problem as a shortest path problem in the same diagram with updated length vectors, as discussed in Propositions~\ref{prop:pathLength} and \ref{prop:pathLengthCompressed}.

\subsection{Results in offline settings}\label{sec:comp_resOffline}

We first presented aggregated results across all combinations of parameters, then discuss how each method is impacted by the choices of hyperparameters. 

\subsubsection{Aggregated results}

\begin{figure}[!h]
\begin{center}
\includegraphics[width=0.45\columnwidth, trim={10.8cm 5.8cm 10.8cm 5.8cm},clip]{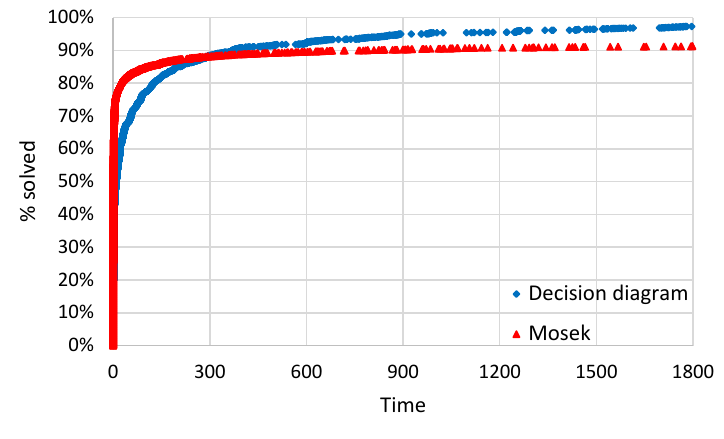}
\caption{Percentage of instances solved as a function of time. The graph summarizes results across 5,400 instances with all combinations of parameters. Mosek solves more instances within a few seconds, but struggles with harder instances, while decision diagrams methods are more effective in the harder instances.}
\label{fig:performanceProfile}
\end{center}
\end{figure}

Table~\ref{tab:resultsFull} shows aggregated results comparing the two different methods in an offline setting. Specifically, the table shows for different values of dimension $n$ and contiguity parameter $\tau$ the metrics of interests. Each row represents an average over five different signals and parameters $k\in \{2,3\}$, $\lambda\in \{0.25,0.50,1.00,2.00,5.00\}$ and $\mu\in \{0.001,0.005,0.010,0.020,0.050,0.100\}$; thus, each row is an average over 300 different instances. The results are summarized in the performance profiles shown in Figure~\ref{fig:performanceProfile}.

\begin{table}[!h]
\begin{center}
\footnotesize
\caption{Computational results as a function of the number of variables $n$ and minimum number of consecutive nonzeros $\tau$. Each row represents an average over 5 different signals $\bm{y}$, $k\in \{2,3\}$, $\lambda\in \{0.25,0.50,1.0,2.00,5.00\}$ and $\mu\in \{0.001,0.005,0.010,0.020,0.050,0.100\}$. Each entry shows the average ``$\pm$" the standard deviation across all instances tested. All times are in seconds.}
\label{tab:resultsFull}
\setlength{\tabcolsep}{1pt}
\begin{tabular}{c c |c c c |c c c c|c}
\hline
 \multirow{2}{*}{$\bm{\tau}$}&\multirow{2}{*}{$\bm{n}$}&\multicolumn{3}{c|}{\textbf{\underline{Mosek}}} & \multicolumn{4}{c|}{\textbf{\underline{Decision diagram}}}& \multicolumn{1}{c}{\textbf{\underline{SOCP}}}\\
&&\textbf{time}&\textbf{nodes}&\textbf{\%}&$\bm{|A|}$&\textbf{time\_dd}&\textbf{time\_sp}&\textbf{\%}&\textbf{time\_hull}\\
 \hline
 \multirow{6}{*}{0}&25&0$\pm$0&38$\pm$103&100\%&15,057$\pm$18,895&1$\pm$0&0$\pm$0&100\%&0$\pm$1\\
&50&1$\pm$4&400$\pm$1,626&100\%&42,052$\pm$57,233&2$\pm$3&0$\pm$0&100\%&3$\pm$7\\
&100&136$\pm$413&18,227$\pm$58,041&96\%&96,032$\pm$133,949&24$\pm$38&0$\pm$0&100\%&26$\pm$67\\
&200&344$\pm$661&8,662$\pm$17,537&85\%&203,992$\pm$287,397&212$\pm$324&0$\pm$0&100\%&78$\pm$191\\
&300&399$\pm$709&4,332$\pm$8,602&82\%&311,952$\pm$440,849&461$\pm$603 &0$\pm$0&93\%&138$\pm$362\\
&500&502$\pm$773&1,188$\pm$2,081&75\%&119,549$\pm$77,521&601$\pm$422&0$\pm$0&65\%&7$\pm$7\\
&&&&&&&&&\\
 \multirow{6}{*}{$5$}&25&0$\pm$0&31$\pm$46&100\%&1,753$\pm$1,524&0$\pm$0&0$\pm$0&100\%&0$\pm$0\\
 &50&1$\pm$2&206$\pm$467&100\%&5,543$\pm$6,068&1$\pm$0&0$\pm$0&100\%&0$\pm$0\\ &100&89$\pm$334&3,379$\pm$12,044&97\%&13,118$\pm$15,164&4$\pm$5&0$\pm$0&100\%&1$\pm$1\\
&200&279$\pm$594&3,016$\pm$6,534&88\%&28,268$\pm$33,360&49$\pm$61&0$\pm$0&100\%&2$\pm$3\\
&300&340$\pm$664&2,595$\pm$5,922&84\%&43,418$\pm$51,557&87$\pm$113&0$\pm$0&100\%&2$\pm$4\\
&500&365$\pm$687&826$\pm$1,746&82\%&73,718$\pm$87,951&459$\pm$508&0$\pm$0&94\%&3$\pm$6\\
&&&&&&&&&\\
 \multirow{6}{*}{$10$}&25&0$\pm$0&14$\pm$20&100\%&622$\pm$312&0$\pm$0&0$\pm$0&100\%&0$\pm$0\\
&50&0$\pm$1&103$\pm$222&100\%&2,195$\pm$1,903&1$\pm$0&0$\pm$0&100\%&0$\pm$0\\ &100&31$\pm$165&1,590$\pm$9,160&100\%&5,360$\pm$5,150&2$\pm$1&0$\pm$0&100\%&0$\pm$0\\
&200&276$\pm$600&2,616$\pm$5,985&88\%&11,690$\pm$11,645&23$\pm$22&0$\pm$0&100\%&1$\pm$1\\
&300&325$\pm$658&1,846$\pm$4,320&84\%&18,020$\pm$18,140&44$\pm$50&0$\pm$0&100\%&1$\pm$1\\
&500&331$\pm$666&468$\pm$1,110&84\%&30,680$\pm$31,131&259$\pm$287&0$\pm$0&100\%&1$\pm$2
\end{tabular}
\end{center}
\end{table}

Observe that the main cost associated with the decision diagram approach is the construction of the diagram itself, while the cost of solving the shortest path problem is negligible. Moreover, the decision diagram approach is competitive with Mosek in all settings, but is particularly effective in constrained instances with large values of $\tau$. Overall, while Mosek seems to be more effective at solving some instances fast (these instances have large values of sparsity parameter $\mu$, see \S\ref{sec:resultsSparsity} for an in-depth analysis), using decision diagrams leads to more instances solved within the 1,800 time limit, and are faster in instances that require more than 300 seconds to solve. Finally, we point that the SOCP reformulation can for the most part be solved within a few seconds, although computational times increase to minutes when the decision diagrams have over 200,000 arcs (see case $\tau=0$, $n\in \{200,300\}$ in Table~\ref{tab:resultsFull}).

\subsubsection{Effect of the sparsity parameter}\label{sec:resultsSparsity} 
\begin{figure}[!h]
\begin{center}
\includegraphics[width=0.5\columnwidth, trim={14cm 6cm 8cm 6cm},clip]{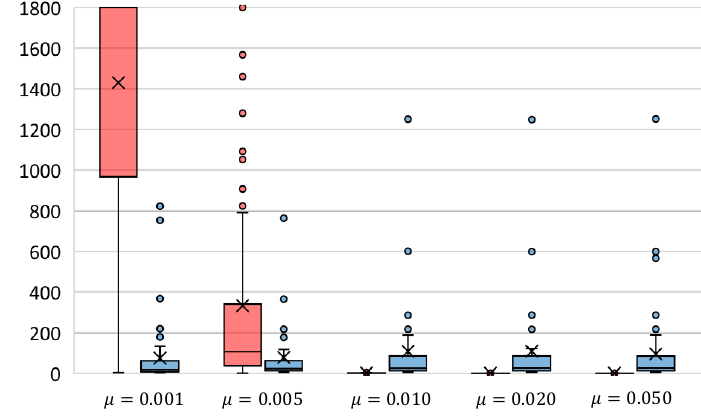}
\caption{Distribution of runtimes of Mosek (red) vs Decision diagram (blue) for $n=200$ as a function of the sparsity parameter $\mu$. Each boxplot represents an average over 5 different signals $\bm{y}$ with $n=200$, $k\in \{2,3\}$ and $\lambda\in \{0.25,0.50,1.0,2.00,5.00\}$.} 
\label{fig:comparisonSparsity}
\end{center}
\end{figure}

Figure~\ref{fig:comparisonSparsity} shows, for instances with $n=200$, the distribution of runtimes of both methods as a function of the sparsity parameter. We observe that the performance of Mosek depends critically on the sparsity parameter $\mu$. When this parameter is large, a branch-and-bound algorithm can easily identify that dense solutions are suboptimal and prune the corresponding branch-and-bound nodes in the tree. In contrast, for small values of this parameter, there are more and denser solutions that achieve good objective values and Mosek struggles to identify an optimal one.

On the other hand, parameter $\mu$, which affects the objective linearly, plays no role in the construction of the decision diagrams. Indeed, as shown in Figure~\ref{fig:comparisonSparsity}, the time to solve the problems via decision diagrams is approximately the same for all values of the sparsity parameter. As a consequence, using our current implementation of the decision diagrams, it is not possible to identify and quickly solve easy instances due to large values of the sparsity parameter. Nonetheless, our implementation of decision diagrams is an order-of-magnitude faster than Mosek in the most challenging instances with $\mu=0.001$.

\subsubsection{Effect of other hyperparameters and cross-validation}\label{sec:resultsOther} Table~\ref{tab:resultsHyperparam} shows, for instances with $n=200$ and for different values of parameters $\tau$, $k$ and $\lambda$, the performance of the different methods. Each row corresponds to the average over five different signals $\bm{y}$ and six different values of sparsity parameter $\mu$. In particular, the computations corresponding to each row mimic a cross-validation procedure to select the best sparsity parameter.

\begin{table}[!h]
\begin{center}
\footnotesize
\caption{Breakdown for $n=200$, as a function of the width $k$, smoothness parameter $\lambda$ and contiguity parameter $\tau$. Each row represents an average over 5 different signals $\bm{y}$ and $\mu\in \{0.001,0.005,0.010,0.020,0.050,0.100\}$. The decision diagram method is able to solve all instances to optimality. The decision diagram constructed for all 30 instances averaged in each row is identical: using decision diagrams, \emph{the 30 instances could be solved by constructing the decision diagram once (incurring a one-time setup cost of time\_dd) and solving the shortest path problem 30 times (with negligible variable cost time\_sp)}.}
\label{tab:resultsHyperparam}
\setlength{\tabcolsep}{2pt}
\begin{tabular}{c c c |c c c |c c c |c}
\hline
\multirow{2}{*}{$\bm{\tau}$}& \multirow{2}{*}{$\bm{k}$}&\multirow{2}{*}{$\bm{\lambda}$}&\multicolumn{3}{c|}{\textbf{\underline{Mosek}}} & \multicolumn{3}{c|}{\textbf{\underline{Decision diagram}}}& \multicolumn{1}{c}{\textbf{\underline{SOCP}}}\\
&&&\textbf{time}&\textbf{nodes}&\textbf{\%}&$\bm{|{A}|}$&\textbf{time\_dd}&\textbf{time\_sp}&\textbf{time\_hull}\\
 \hline
\multirow{11}{*}{0} &\multirow{5}{*}{2}&0.25&47$\pm$130&1,119$\pm$3,004&100\%&10,965$\pm$0&7$\pm$1&0$\pm$0&0$\pm$0\\
&&0.5&421$\pm$741&7,628$\pm$13,572&83\%&16,749$\pm$0&12$\pm$1&0$\pm$0&0$\pm$0\\
 &&1.0&457$\pm$759&15,250$\pm$26,760&77\%&30,963$\pm$0&23$\pm$4&0$\pm$0&1$\pm$0\\
 &&2.0&474$\pm$769&11,672$\pm$19,010&77\%&51,923$\pm$0&44$\pm$6&0$\pm$0&2$\pm$0\\
 &&5.0&351$\pm$677&8,571$\pm$16,533&83\%&88,491$\pm$0&79$\pm$14&0$\pm$0&4$\pm$1\\
 &&&&&&&\\
 &\multirow{5}{*}{3}&0.25&18$\pm$45&373$\pm$848&100\%&56,789$\pm$0&55$\pm$9&0$\pm$0&2$\pm$1\\
 &&0.5&280$\pm$563&5,920$\pm$12,841&93\%&107,591$\pm$0&104$\pm$23&0$\pm$0&6$\pm$2\\
 &&1.0&511$\pm$780&14,590$\pm$23,573&77\%&233,917$\pm$0&248$\pm$58&0$\pm$0&29$\pm$6\\
 &&2.0&471$\pm$762&11,448$\pm$18,690&77\%&478,889$\pm$0&524$\pm$112&0$\pm$0&125$\pm$26\\
 &&5.0&408$\pm$720&10,048$\pm$17,800&83\%&963,643$\pm$0&1,026$\pm$268&0$\pm$0&611$\pm$190\\
 &&&&&&&&&\\
 \hline
 
 \multirow{11}{*}{5} &\multirow{5}{*}{2}&0.25&136$\pm$298&1,434$\pm$3,205&100\%&3,124$\pm$0&5$\pm$1&0$\pm$0&0$\pm$0\\
 &&0.5&367$\pm$681&3,675$\pm$6,815&83\%&5,420$\pm$0&8$\pm$2&0$\pm$0&0$\pm$0\\
 &&1.0&364$\pm$680&3,505$\pm$6,543&83\%&8,440$\pm$0&12$\pm$2&0$\pm$0&0$\pm$0\\
 &&2.0&319$\pm$675&3,571$\pm$8,169&83\%&13,136$\pm$0&19$\pm$2&0$\pm$0&1$\pm$0\\
 &&5.0&304$\pm$681&3,779$\pm$8,687&83\%&23,141$\pm$0&35$\pm$3&0$\pm$0&1$\pm$0\\
 &&&&&&&&&\\
 &\multirow{5}{*}{3}&0.25&92$\pm$200&946$\pm$2,071&100\%&11,289$\pm$0&20$\pm$3&0$\pm$0&0$\pm$0\\
 &&0.5&194$\pm$431&2,789$\pm$5,331&97\%&15,442$\pm$0&26$\pm$6&0$\pm$0&1$\pm$0\\
 &&1.0&386$\pm$685&3,703$\pm$6,576&83\%&30,520$\pm$0&54$\pm$4&0$\pm$0&1$\pm$0\\
 &&2.0&327$\pm$674&4,261$\pm$8,957&83\%&53,325$\pm$0&93$\pm$11&0$\pm$0&3$\pm$1\\
 &&5.0&304$\pm$681&2,499$\pm$5,591&83\%&118,842$\pm$0&214$\pm$10&0$\pm$0&12$\pm$1\\
 &&&&&&&&&\\
 \hline
 
 \multirow{11}{*}{10} &\multirow{5}{*}{2}&0.25&125$\pm$230&1,355$\pm$2,522&100\%&3,661$\pm$0&7$\pm$1&0$\pm$0&0$\pm$0\\
 &&0.5&331$\pm$672&3,497$\pm$7,506&83\%&3,846$\pm$0&8$\pm$1&0$\pm$0&0$\pm$0\\
 &&1.0&325$\pm$673&2,775$\pm$5,786&83\%&4,766$\pm$0&10$\pm$2&0$\pm$0&0$\pm$0\\
 &&2.0&289$\pm$636&3,897$\pm$9,178&87\%&6,431$\pm$0&12$\pm$2&0$\pm$0&1$\pm$0\\
 &&5.0&303$\pm$681&2,081$\pm$4,809&83\%&11,903$\pm$0&21$\pm$1&0$\pm$0&1$\pm$0\\
 &&&&&&&\\
 &\multirow{5}{*}{3}&0.25&145$\pm$330&1,458$\pm$3,223&100\%&5,693$\pm$0&13$\pm$3&0$\pm$0&0$\pm$0\\
&&0.5&348$\pm$672&2,908$\pm$5,827&83\%&6,439$\pm$0&15$\pm$3&0$\pm$0&0$\pm$0\\
 &&1.0&322$\pm$674&2,712$\pm$5,623&83\%&10,310$\pm$0&23$\pm$3&0$\pm$0&1$\pm$0\\
 &&2.0&287$\pm$630&3,296$\pm$7,513&90\%&20,398$\pm$0&39$\pm$3&0$\pm$0&1$\pm$0\\
 &&5.0&285$\pm$648&2,181$\pm$5,171&87\%&43,448$\pm$0&84$\pm$3&0$\pm$0&2$\pm$0
\end{tabular}
\end{center}
\end{table}

Since the averages in each row are taken over parameters concerning the linear coefficients of the objective, which do not influence the construction of the decision diagram, the ensuing diagram is the same for all instances. In particular, during cross-validation, the decision diagram needs to be constructed only once and then can be reused with minimal cost. In contrast, to solve these instances with Mosek, repeated calls to a branch-and-bound solver need to be done, with limited potential for reoptimization. We also observe that, as expected, decision diagrams are substantially more effective for smaller values of the regularization parameter $\lambda$ (leading to better conditioned matrices) and width parameter $k$. As observed previously, decision diagrams are also more effective for heavily constrained instances with larger values of parameter $\tau$.  

\begin{figure}[!h]
    \centering
    \subfigure[Unconstrained instances: $\tau=0$]{
    \includegraphics[width=0.35\textwidth,trim={11.5cm 5cm 11.5cm 5cm},clip]{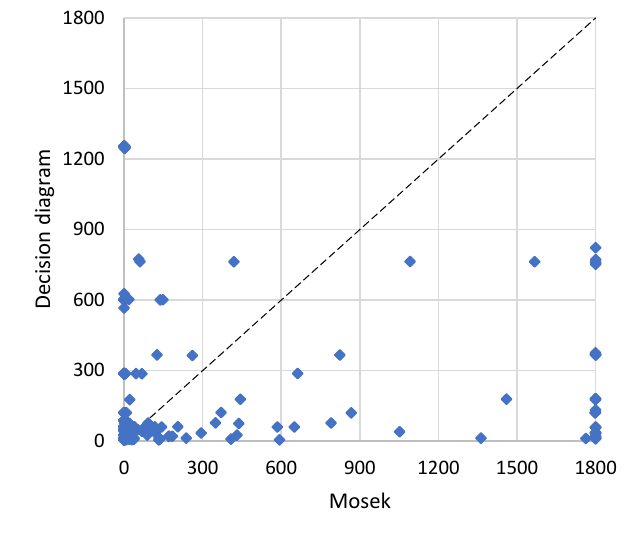}
    }
    \subfigure[Constrained Instances with $\tau\in \{5,10\}$]{
    \includegraphics[width=0.35\textwidth,trim={11.5cm 5cm 11.5cm 5cm},clip]{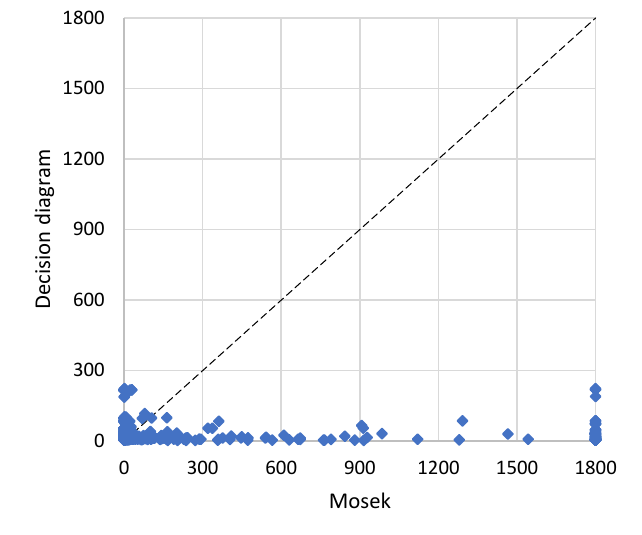}
    }
    
    \caption{Plots comparison solution times of all instances with $n=200$. Each dot corresponds to an instance, with the horizontal coordinate representing the time used by Mosek and the vertical coordinate corresponding to the time used by the decision diagram method.}
    \label{fig:timesComparison}
\end{figure}

One of the main advantages of decision diagram approaches is the capabilities of reoptimization, as shown in Table~\ref{tab:resultsHyperparam} and, later, with computations in the online setting in \S\ref{sec:comp_resOnline}. Nonetheless, we emphasize that they can improve upon off-the-shelf solvers even in the context of a single instance. Figure~\ref{fig:timesComparison} depicts the time required to solve each individual instance with $n=200$ from scratch using either Mosek or via decision diagrams. We see that even in unconstrained instances, where decision diagrams are less effective, the method is competitive with Mosek: Indeed, while Mosek is able to solve some instances (with large values of $\mu$) almost instantly, it also hits time limits in several instances. Decision diagrams, on the other hand, are able to solve all instances consistently, and are competitive with the off-the-shelf solver. In constrained instances with $\tau\in \{5,10\}$, using decision diagrams clearly results in better performance, solving all instances to optimality within 300 seconds, and faster than Mosek in most cases.

\subsection{Results in online settings}\label{sec:comp_resOnline}

Table~\ref{tab:resultsOnline} shows results for the online setting, using the same hyperparameters combinations as in the offline setting. The table shows for each combination of parameters the number of arcs of the decision diagram \underline{$|A|$} and the setup time \underline{time\_dd} (corresponding exactly to the offline setting), the average time per instance \underline{time\_sp} and the total time \underline{time\_total} required to solve all 6,823 instances online (not including setup time). All times reported are in seconds. The results are also summarized in Figure~\ref{fig:timesOnline} in the introduction.

\begin{table}[H]
\begin{center}
\footnotesize
\caption{Online instances, each one requiring the sequential solution of 6,823 MIOs \eqref{eq:denoising} with $n=200$ (corresponding, for each point, to the most recent $200$ observations).}
\label{tab:resultsOnline}
\setlength{\tabcolsep}{2pt}
\begin{tabular}{ c c |c c | c c }
\hline
 \multirow{2}{*}{$\bm{k}$}&\multirow{2}{*}{$\bm{\lambda}$}& \multicolumn{2}{c|}{\textbf{\underline{Setup time}}}& \multicolumn{2}{c}{\textbf{\underline{Online time}}}\\
&&$\bm{|A|}$&\textbf{time\_dd (s)}&\textbf{time\_sp(s)}&\textbf{time\_total(s)}\\
 \hline
\multirow{5}{*}{2}&0.25&10,965&7&0.001&7\\
&0.5&16,749&11&0.002&11\\
&1.0&30,963&24&0.004&30\\
&2.0&51,923&32&0.006&43\\
&5.0&88,491&62&0.013&88\\
 &&&&&\\
 \multirow{5}{*}{3}&0.25&56,789&40&0.007&48\\
&0.5&107,591&81&0.016&107\\
&1.0&233,917&184&0.035&239\\
&2.0&478,889&409&0.079&539\\
&5.0&963,643&864&0.185&1,261
\end{tabular}
\end{center}
\end{table}

Decision diagrams are able to solve the problems in milliseconds. In instances with $k=2$ and $\lambda$ small, the solution times per instance can be as small as one millisecond per instance. Even in instances with larger $k$ and $\lambda$, the solution times are under $0.2$ seconds on average and under $0.4$ seconds in the worst case.

\section*{Acknowledgmenets}Dr. Lozano gratefully acknowledges the support of the \emph{Air Force Office of Scientific Research} under grant FA9550-22-1-0236. Andr\'es G\'omez is supported, in part, by grant FA9550-22-1-0369 from the Air Force Office of Scientific Research and grant 2152777 from the National Science Foundation.


\bibliographystyle{informs2014}
\bibliography{references.bib}

\begin{thebibliography}{46}
\providecommand{\natexlab}[1]{#1}
\providecommand{\url}[1]{\texttt{#1}}
\providecommand{\urlprefix}{URL }

\bibitem[{Atamt{\"u}rk et~al.(2021)Atamt{\"u}rk, G{\'o}mez, \protect\BIBand{}
  Han}]{atamturk2021sparse}
Atamt{\"u}rk A, G{\'o}mez A, Han S (2021) Sparse and smooth signal estimation:
  Convexification of $\ell_0$-formulations. \emph{The Journal of Machine
  Learning Research} 22(1):2370--2412.

\bibitem[{Bergman et~al.(2016{\natexlab{a}})Bergman, Cire, Van~Hoeve,
  \protect\BIBand{} Hooker}]{bergman2016}
Bergman D, Cire AA, Van~Hoeve WJ, Hooker J (2016{\natexlab{a}}) \emph{Decision
  diagrams for optimization}, volume~1 (Springer).

\bibitem[{Bergman et~al.(2014)Bergman, Cire, van Hoeve, \protect\BIBand{}
  Hooker}]{BerCirHoeHoo14}
Bergman D, Cire AA, van Hoeve WJ, Hooker JN (2014) Optimization bounds from
  binary decision diagrams. \emph{INFORMS Journal on Computing} 26(2):253--268.

\bibitem[{Bergman et~al.(2016{\natexlab{b}})Bergman, Cire, van Hoeve,
  \protect\BIBand{} Hooker}]{BerCirHoeHoo16}
Bergman D, Cire AA, van Hoeve WJ, Hooker JN (2016{\natexlab{b}}) Discrete
  optimization with decision diagrams. \emph{INFORMS Journal on Computing}
  28(1):47--66.

\bibitem[{Bertsimas et~al.(2024)Bertsimas, Digalakis~Jr, Li, \protect\BIBand{}
  Lami}]{bertsimas2024slowly}
Bertsimas D, Digalakis~Jr V, Li ML, Lami OS (2024) Slowly varying regression
  under sparsity. \emph{Operations Research} .

\bibitem[{Bertsimas et~al.(2016)Bertsimas, King, \protect\BIBand{}
  Mazumder}]{bertsimas2016best}
Bertsimas D, King A, Mazumder R (2016) Best subset selection via a modern
  optimization lens. \emph{The Annals of Statistics} 813--852.

\bibitem[{Bhathena et~al.(2024)Bhathena, Fattahi, G{\'o}mez, \protect\BIBand{}
  K{\"u}{\c{c}}{\"u}kyavuz}]{bhathena2024parametric}
Bhathena A, Fattahi S, G{\'o}mez A, K{\"u}{\c{c}}{\"u}kyavuz S (2024) A
  parametric approach for solving convex quadratic optimization with indicators
  over trees. \emph{arXiv preprint arXiv:2404.08178} .

\bibitem[{Campajola et~al.(2022)Campajola, Gangi, Lillo, \protect\BIBand{}
  Tantari}]{campajola2022modelling}
Campajola C, Gangi DD, Lillo F, Tantari D (2022) Modelling time-varying
  interactions in complex systems: the score driven kinetic ising model.
  \emph{Scientific Reports} 12(1):19339.

\bibitem[{Candes et~al.(2008)Candes, Wakin, \protect\BIBand{}
  Boyd}]{candes2008enhancing}
Candes EJ, Wakin MB, Boyd SP (2008) Enhancing sparsity by reweighted $\ell_1$
  minimization. \emph{Journal of Fourier Analysis and Applications}
  14:877--905.

\bibitem[{Castro et~al.(2022)Castro, Cire, \protect\BIBand{}
  Beck}]{castro2022decision}
Castro MP, Cire AA, Beck JC (2022) Decision diagrams for discrete optimization:
  A survey of recent advances. \emph{INFORMS Journal on Computing}
  34(4):2271--2295.

\bibitem[{Chen et~al.(2001)Chen, Donoho, \protect\BIBand{}
  Saunders}]{chen2001atomic}
Chen SS, Donoho DL, Saunders MA (2001) Atomic decomposition by basis pursuit.
  \emph{SIAM review} 43(1):129--159.

\bibitem[{Cire \protect\BIBand{} van Hoeve(2013)}]{CirHoe13}
Cire AA, van Hoeve WJ (2013) Multivalued decision diagrams for sequencing
  problems. \emph{Operations Research} 61(6):1411--1428.

\bibitem[{Das \protect\BIBand{} Kempe(2008)}]{das2008algorithms}
Das A, Kempe D (2008) Algorithms for subset selection in linear regression.
  \emph{Proceedings of the Fortieth Annual ACM Symposium on Theory of
  Computing}, 45--54.

\bibitem[{Davarnia \protect\BIBand{} Van~Hoeve(2021)}]{davarnia2020}
Davarnia D, Van~Hoeve WJ (2021) Outer approximation for integer nonlinear
  programs via decision diagrams. \emph{Mathematical Programming} 187:111--150.

\bibitem[{Demko et~al.(1984)Demko, Moss, \protect\BIBand{}
  Smith}]{demko1984decay}
Demko S, Moss WF, Smith PW (1984) Decay rates for inverses of band matrices.
  \emph{Mathematics of Computation} 43(168):491--499.

\bibitem[{Donoho et~al.(2005)Donoho, Elad, \protect\BIBand{}
  Temlyakov}]{donoho2005stable}
Donoho DL, Elad M, Temlyakov VN (2005) Stable recovery of sparse overcomplete
  representations in the presence of noise. \emph{IEEE Transactions on
  Information Theory} 52(1):6--18.

\bibitem[{Dunn et~al.(2018)Dunn, Runge, \protect\BIBand{}
  Snyder}]{dunn2018wearables}
Dunn J, Runge R, Snyder M (2018) Wearables and the medical revolution.
  \emph{Personalized Medicine} 15(5):429--448.

\bibitem[{Frangioni \protect\BIBand{} Gentile(2006)}]{frangioni2006perspective}
Frangioni A, Gentile C (2006) Perspective cuts for a class of convex 0--1 mixed
  integer programs. \emph{Mathematical Programming} 106:225--236.

\bibitem[{Friedrich et~al.(2017)Friedrich, Zhou, \protect\BIBand{}
  Paninski}]{friedrich2017fast}
Friedrich J, Zhou P, Paninski L (2017) Fast online deconvolution of calcium
  imaging data. \emph{PLoS Computational Biology} 13(3):e1005423.

\bibitem[{G{\'o}mez(2021)}]{gomez2021outlier}
G{\'o}mez A (2021) Outlier detection in time series via mixed-integer conic
  quadratic optimization. \emph{SIAM Journal on Optimization} 31(3):1897--1925.

\bibitem[{G{\"u}nl{\"u}k \protect\BIBand{}
  Linderoth(2010)}]{gunluk2010perspective}
G{\"u}nl{\"u}k O, Linderoth J (2010) Perspective reformulations of mixed
  integer nonlinear programs with indicator variables. \emph{Mathematical
  Programming} 124:183--205.

\bibitem[{Guo et~al.(2016)Guo, Hu, Jing, \protect\BIBand{}
  Zhang}]{guo2016spline}
Guo J, Hu J, Jing BY, Zhang Z (2016) Spline-lasso in high-dimensional linear
  regression. \emph{Journal of the American Statistical Association}
  111(513):288--297.

\bibitem[{Han \protect\BIBand{} G{\'o}mez(2021)}]{han2021compact}
Han S, G{\'o}mez A (2021) Compact extended formulations for low-rank functions
  with indicator variables. \emph{arXiv preprint arXiv:2110.14884} .

\bibitem[{Han et~al.(2022)Han, G{\'o}mez, \protect\BIBand{}
  Pang}]{han2022polynomial}
Han S, G{\'o}mez A, Pang JS (2022) On polynomial-time solvability of
  combinatorial markov random fields. \emph{arXiv preprint arXiv:2209.13161} .

\bibitem[{Hazimeh et~al.(2022)Hazimeh, Mazumder, \protect\BIBand{}
  Saab}]{hazimeh2022sparse}
Hazimeh H, Mazumder R, Saab A (2022) Sparse regression at scale:
  Branch-and-bound rooted in first-order optimization. \emph{Mathematical
  Programming} 196(1-2):347--388.

\bibitem[{Hodrick \protect\BIBand{} Prescott(1997)}]{hodrick1997postwar}
Hodrick RJ, Prescott EC (1997) Postwar us business cycles: an empirical
  investigation. \emph{Journal of Money, Credit, and Banking} 1--16.

\bibitem[{Hoerl \protect\BIBand{} Kennard(1970)}]{hoerl1970ridge}
Hoerl AE, Kennard RW (1970) Ridge regression: Biased estimation for
  nonorthogonal problems. \emph{Technometrics} 12(1):55--67.

\bibitem[{Jewell \protect\BIBand{} Witten(2018)}]{jewell2018exact}
Jewell S, Witten D (2018) Exact spike train inference via $\ell_0$
  optimization. \emph{The Annals of Applied Statistics} 12(4):2457.

\bibitem[{Kim et~al.(2009)Kim, Koh, Boyd, \protect\BIBand{}
  Gorinevsky}]{kim2009ell_1}
Kim SJ, Koh K, Boyd S, Gorinevsky D (2009) $\ell\_1$ trend filtering.
  \emph{SIAM Review} 51(2):339--360.

\bibitem[{Lin et~al.(2014)Lin, Pham, \protect\BIBand{}
  Ruszczy{\'n}ski}]{lin2014alternating}
Lin X, Pham M, Ruszczy{\'n}ski A (2014) Alternating linearization for
  structured regularization problems. \emph{The Journal of Machine Learning
  Research} 15(1):3447--3481.

\bibitem[{Liu et~al.(2023)Liu, Fattahi, G{\'o}mez, \protect\BIBand{}
  K{\"u}{\c{c}}{\"u}kyavuz}]{liu2023graph}
Liu P, Fattahi S, G{\'o}mez A, K{\"u}{\c{c}}{\"u}kyavuz S (2023) A graph-based
  decomposition method for convex quadratic optimization with indicators.
  \emph{Mathematical Programming} 200(2):669--701.

\bibitem[{Lozano \protect\BIBand{} Smith(2022)}]{lozano2018Binary}
Lozano L, Smith JC (2022) A binary decision diagram based algorithm for solving
  a class of binary two-stage stochastic programs. \emph{Mathematical
  Programming} 191:381--404.

\bibitem[{Lu \protect\BIBand{} Shiou(2002)}]{lu2002inverses}
Lu TT, Shiou SH (2002) Inverses of 2$\times$ 2 block matrices. \emph{Computers
  \& Mathematics with Applications} 43(1-2):119--129.

\bibitem[{MacNeil \protect\BIBand{} Bodur(2024)}]{macneil2024leveraging}
MacNeil M, Bodur M (2024) Leveraging decision diagrams to solve two-stage
  stochastic programs with binary recourse and logical linking constraints.
  \emph{European Journal of Operational Research} 315(1):228--241.

\bibitem[{Mammen \protect\BIBand{} Van De~Geer(1997)}]{mammen1997locally}
Mammen E, Van De~Geer S (1997) Locally adaptive regression splines. \emph{The
  Annals of Statistics} 25(1):387--413.

\bibitem[{Mazumder et~al.(2023)Mazumder, Radchenko, \protect\BIBand{}
  Dedieu}]{mazumder2023subset}
Mazumder R, Radchenko P, Dedieu A (2023) Subset selection with shrinkage:
  Sparse linear modeling when the snr is low. \emph{Operations Research}
  71(1):129--147.

\bibitem[{Nesterov \protect\BIBand{} Nemirovskii(1994)}]{nesterov1994interior}
Nesterov Y, Nemirovskii A (1994) \emph{Interior-point polynomial algorithms in
  convex programming} (SIAM).

\bibitem[{Rinaldo(2009)}]{rinaldo2009properties}
Rinaldo A (2009) Properties and refinements of the fused lasso. \emph{The
  Annals of Statistics} 37:2922--2952.

\bibitem[{Ruppert(2002)}]{ruppert2002selecting}
Ruppert D (2002) Selecting the number of knots for penalized splines.
  \emph{Journal of Computational and Graphical Statistics} 11(4):735--757.

\bibitem[{Tibshirani et~al.(2005)Tibshirani, Saunders, Rosset, Zhu,
  \protect\BIBand{} Knight}]{tibshirani2005sparsity}
Tibshirani R, Saunders M, Rosset S, Zhu J, Knight K (2005) Sparsity and
  smoothness via the fused lasso. \emph{Journal of the Royal Statistical
  Society Series B: Statistical Methodology} 67(1):91--108.

\bibitem[{Vogelstein et~al.(2010)Vogelstein, Packer, Machado, Sippy, Babadi,
  Yuste, \protect\BIBand{} Paninski}]{vogelstein2010fast}
Vogelstein JT, Packer AM, Machado TA, Sippy T, Babadi B, Yuste R, Paninski L
  (2010) Fast nonnegative deconvolution for spike train inference from
  population calcium imaging. \emph{Journal of Neurophysiology}
  104(6):3691--3704.

\bibitem[{Wei et~al.(2024)Wei, Atamt{\"u}rk, G{\'o}mez, \protect\BIBand{}
  K{\"u}{\c{c}}{\"u}kyavuz}]{wei2023convex}
Wei L, Atamt{\"u}rk A, G{\'o}mez A, K{\"u}{\c{c}}{\"u}kyavuz S (2024) On the
  convex hull of convex quadratic optimization problems with indicators.
  \emph{Mathematical Programming} 204(1):703--737.

\bibitem[{Xie \protect\BIBand{} Deng(2020)}]{xie2020scalable}
Xie W, Deng X (2020) Scalable algorithms for the sparse ridge regression.
  \emph{SIAM Journal on Optimization} 30(4):3359--3386.

\bibitem[{Yan et~al.(2014)Yan, Paynabar, \protect\BIBand{} Shi}]{yan2014image}
Yan H, Paynabar K, Shi J (2014) Image-based process monitoring using low-rank
  tensor decomposition. \emph{IEEE Transactions on Automation Science and
  Engineering} 12(1):216--227.

\bibitem[{Yan et~al.(2017)Yan, Paynabar, \protect\BIBand{}
  Shi}]{yan2017anomaly}
Yan H, Paynabar K, Shi J (2017) Anomaly detection in images with smooth
  background via smooth-sparse decomposition. \emph{Technometrics}
  59(1):102--114.

\bibitem[{Zou \protect\BIBand{} Qiu(2009)}]{zou2009multivariate}
Zou C, Qiu P (2009) Multivariate statistical process control using lasso.
  \emph{Journal of the American Statistical Association} 104(488):1586--1596.

\end{thebibliography}

\clearpage
\appendix

\section{Proofs}\label{sec:proofs}

\proof{(Proposition~\ref{prop:correctReduced})}
We prove the result by induction. To simplify the notation, for any depth $\ell$ and indexes $i,j\in [n]$ we let $\bm{(Q^\ell)^\dagger}=\left(\bm{Q}\circ\bm{\bar z^\ell}(\bm{\bar z^\ell})^\top\right)^\dagger$, we denote by $\bm{\bar W^\ell_i}$ and $\bm{(Q^\ell)_i^\dagger}$ denote the $i$-th columns of $\bm{\bar W^\ell}$ and $\bm{(Q^\ell)^\dagger}$, respectively, and by $(Q^\ell)_{ij}^\dagger$ the $(i,j)$-th entry of $\bm{(Q^\ell)^\dagger}$. 

\noindent \textbf{Base case} If $\ell=1$, $s^\ell=\left\{\bm{0}_{n\times n}\right\}$ and the properties are trivially satisfied.

\noindent \textbf{Inductive step} Assume that the results hold for path $(a_1,a_2,\dots,a_{\ell-1})$, and we prove them for path $(a_1,a_2,\dots,a_{\ell-1},a_{\ell}).$ There are three cases, depending on the depth $\ell$ and the assignment $\hat z_{\ell}$. \newline
$\bullet$ If $\ell\geq \pi_j\Leftrightarrow \ell+1>\pi_j$, then $\bar W_{ij}^{\ell+1}=0$.\newline
$\bullet$ If $\ell<\pi_j\Leftrightarrow \ell+1\leq \pi_j$ and $\hat z_\ell=0$, then $$\bar W_{ij}^{\ell+1}=\bar W_{ij}^\ell=(Q^\ell)_{ij}^\dagger=(Q^{\ell+1})_{ij}^\dagger,$$
where the first equality follows since the state does not change, the second equality follows due to the induction hypothesis, and the third equality follows since $\hat z_\ell=0\implies \bm{\bar z^{\ell}}=\bm{\bar z^{\ell+1}}$. Thus the first property is satisfied. Moreover, in this case $\bm{u}=\bm{0}$, also satisfying the second property.\newline
$\bullet$ If $\ell<\pi_j\Leftrightarrow \ell+1\leq \pi_j$ and $\hat z_\ell=1$, then we find that
\begin{align*}\bm{u}&=\frac{1}{\sqrt{Q_{\ell\ell}-\sum_{p=1}^n\sum_{q=1}^n\bar W_{ pq }^\ell Q_{p \ell}Q_{q \ell}}}\left(-\sum_{q=1}^n \bm{\bar W_q^\ell}Q_{q \ell}+\bm{e_\ell}\right)\\
&=\frac{1}{\sqrt{Q_{\ell\ell}-\sum_{p=1}^n\sum_{q\in [n]:\pi_q\geq \ell}(Q^\ell)_{p q}^\dagger Q_{p \ell}Q_{q \ell}}}\left(-\sum_{q\in [n]:\pi_q\geq \ell} \bm{(Q^\ell)_q^\dagger}Q_{q \ell}+\bm{e_\ell}\right)\tag{$\because $ induction hypothesis}\\
&=\frac{1}{\sqrt{Q_{\ell\ell}-\sum_{p=1}^n\sum_{q=1}^n(Q^\ell)_{p q}^\dagger Q_{p \ell}Q_{q \ell}}}\left(-\sum_{q=1}^n \bm{(Q^\ell)_q^\dagger}Q_{q \ell}+\bm{e_\ell}\right)\tag{$\because Q_{q \ell}=0$ for $\pi_q<\ell$}.\end{align*}
From this point the proof is identical to the corresponding step in Proposition~\ref{prop:state}: from Corollary~\ref{cor:inversion}, we find that $\bm{uu^\top}=\left(\bm{Q}\circ\bm{\bar z^{\ell+1}}(\bm{\bar z^{\ell+1}})^\top\right)^\dagger-\left(\bm{Q}\circ\bm{\bar z^\ell}(\bm{\bar z^\ell})^\top\right)^\dagger$, where $\bm{z^{\ell+1}}=\bm{z_\ell}+\bm{e_\ell}$, thus proving the second property. The first property follows immediately. 
\endproof

\section{On numerical precision}\label{sec:numPrec}

\begin{figure}[!h]
    \centering
    \subfigure[Relative objective gap, avg=-0.0008\%, std=0.0069\%]{
    \includegraphics[width=0.48\textwidth,trim={10cm 6cm 10cm 6cm},clip]{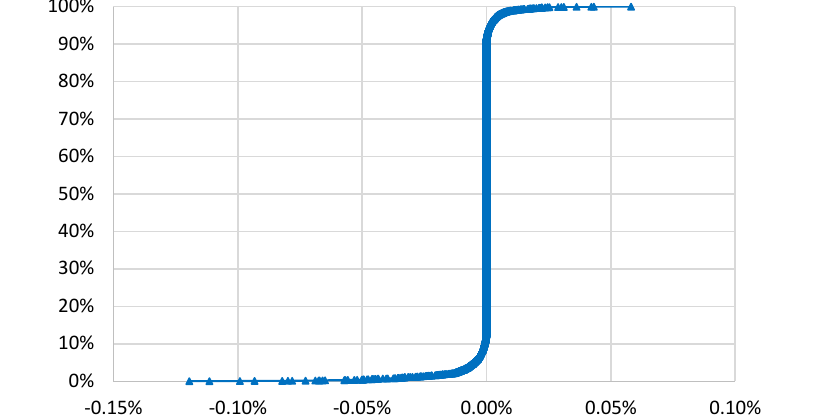}
    }
    \subfigure[Relative solution gap, avg=0.0001\%, std=0.0012\%]{
    \includegraphics[width=0.48\textwidth,trim={10cm 6cm 10cm 6cm},clip]{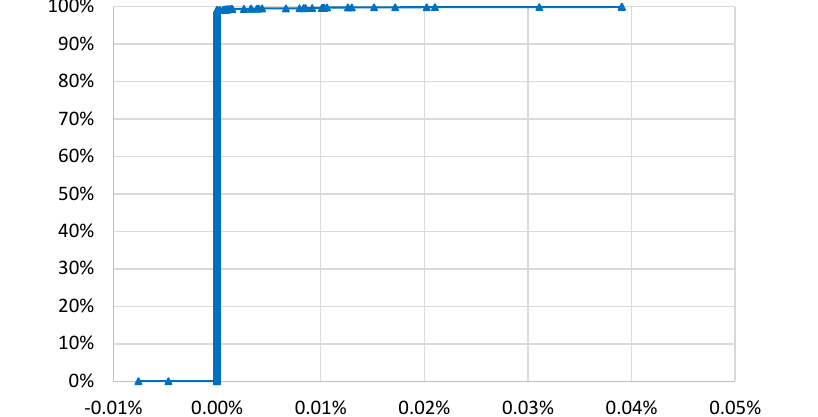}
    }
    
    \caption{Numerical precision comparison between decision diagrams and Mosek mixed-integer solver.}
    \label{fig:DDvsMosekPrecision}
\end{figure}

In this section we comment on the numerical precision of the proposed methods. Note that, in general, methods for MIQO are subject to precision errors, and the proposed approach is no exception. In particular, the $\epsilon$-exact decision diagrams introduce round-off errors due to merging of states that, although close, are not equal. To evaluate the magnitude of the precision errors introduced, we compare the solutions obtained by the proposed method and with those reported by Mosek, and report two metrics. We report the relative objective gap computed as $\frac{\text{obj}_\text{dd}-\text{obj}_\text{msk}}{\text{obj}_\text{msk}}$, where $\text{obj}_\text{dd}$ and $\text{obj}_\text{msk}$ are the objective values reported by the decision diagram and Mosek, respectively. We also report the relative solution gap computed as $$\frac{\bm{c^\top z_{\textbf{dd}}}-\frac{1}{2}\bm{d^\top}\left(\bm{Q}\circ \bm{z_{\textbf{dd}}z_{\textbf{dd}}^\top}\right)^\dagger\bm{d}-\bm{c^\top z_{\textbf{msk}}}+\frac{1}{2}\bm{d^\top}\left(\bm{Q}\circ \bm{z_{\textbf{msk}}z_{\textbf{msk}}^\top}\right)^\dagger\bm{d}}{\bm{c^\top z_{\textbf{msk}}}-\frac{1}{2}\bm{d^\top}\left(\bm{Q}\circ \bm{z_{\textbf{msk}}z_{\textbf{msk}}^\top}\right)^\dagger\bm{d}},$$ where $z_{\textbf{dd}}$ and $z_{\textbf{msk}}$ are the solutions found by the decision diagram and Mosek, respectively. Note that in the second metric we are manually verifying the objective value associated with the solution reported by each solver. Figure~\ref{fig:DDvsMosekPrecision} reports the distribution of these two metrics across the 4,827 instances that both methods solved to optimality in our computational experiments. 

Both figures show the difference of quality between the two approach is minimal. In 99.0\% of the instances the solution gap is zero, that is, the solutions reported by both methods are identical. Moreover, in 99.7\% of the instances the solution gap between the approaches is less than $10^{-4}$ (the default optimality gap tolerance for Mosek): in other words, the solutions produced by decision diagrams could be certified by Mosek as optimal in a vast majority of the instances. We point out that in inference problems such as \eqref{eq:inference} the matrices involved are relatively well-conditioned, see Remark~\ref{rem:condition}. As a consequence our code, even though it does not implement any advanced numerical techniques to control numerical errors, still delivers optimal solutions in almost all cases. 

We now turn our attention to solutions produced by the convex reformulation in Theorem~\ref{theo:hullDD}. Observe that solutions produced by this method could be subject to two additional sources of numerical errors. First, the formulation involves vectors $\bm{u}$ themselves instead of outer products of these vectors, hence the computation of the square root results in unavoidable round-off errors (interestingly, the intermediate SDP reformulation \eqref{eq:hullSdp} does not have this issue). Second, since this solution is then obtained by using a nonlinear solver, the method is affected by the numerical precision of the solver. We point out that the SOCP formulations are large, involving anywhere from a few thousand to close to a million additional variables: even modern codes, that deliver extremely high-quality solutions in problems with hundreds of variables, can struggle to avoid round-off errors with problems of this size. Finally, although not a numerical precision error, we point out that in instances with multiple optimal solutions, interior point methods for SOCPs would tend to deliver one that is in the interior of the convex hull of the set of optimal solutions: this solution might not be integral in variables $\bm{z}$.

\begin{figure}[!h]
    \centering
    \subfigure[Relative obj gap, avg=-0.0601\%, std=0.4020\%]{
    \includegraphics[width=0.44\textwidth,trim={10cm 6cm 10cm 6cm},clip]{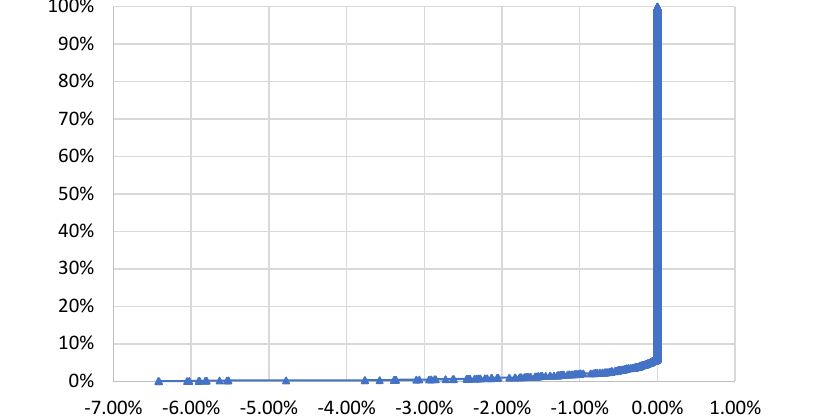}
    }
    \subfigure[Relative solution gap, avg=-0.0516\%, std=0.3742\%]{
    \includegraphics[width=0.44\textwidth,trim={10cm 6cm 10cm 6cm},clip]{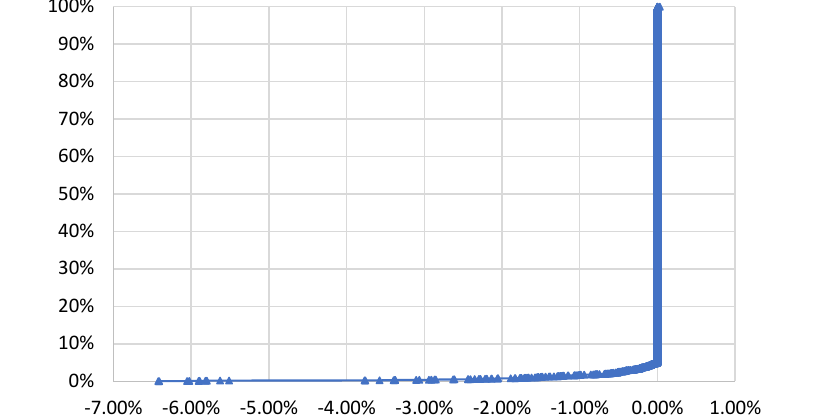}
    }
    
    \caption{Numerical precision comparison between decision diagrams and Mosek solving the convex reformulation induced by Theorem~\ref{theo:hullDD}.}
    \label{fig:DDvsRelPrecision}
\end{figure}

Figure~\ref{fig:DDvsRelPrecision} presents the relative objective and solution gaps resulting from solving \eqref{eq:miqo} via a shortest path on the decision diagram vs solving the convex reformulation ensuing from Theorem~\ref{theo:hullDD}. These metrics are computed in the same fashion, but $\text{obj}_{\text{msk}}$ and $\text{z}_{\textbf{msk}}$ are replaced with the objective and solution obtained from the convex reformulation instead of the mixed-integer optimization solver. In 94.9\% the solution gap is exactly zero, showing that the convex reformulation produces optimal solution in a majority of the instances as well. However, in 1.7\% of the instances, the convex reformulation produces a solution with solution gap smaller than -1\%, thus suboptimal by a noticeable margin. We conclude that using the proposed shortest path method to solve problems with decision diagrams is not only faster, but results in more precise solutions as well. 

\end{document}